\def\fileversion{v2.25a} \def\filedate{2009/10/10}
\documentclass[edeposit,fullpage]{uiucthesis2009}

\usepackage{url} 
\usepackage{graphicx} 
\usepackage{listings} 
\usepackage[usenames,dvipsnames]{color} 

\usepackage{amsmath} 
\usepackage{amsfonts} 
\usepackage{amsthm} 
\usepackage{amssymb} 
\usepackage{amscd} 
\usepackage{cancel} 
\usepackage{sagetex}
\usepackage{boxedminipage}
\usepackage{morefloats}

\newsavebox{\fmbox}

\newcommand{\citep}{\cite} 

\newcommand{\gkpSI}[2]{\ensuremath{\genfrac{\lbrack}{\rbrack}{0pt}{}{#1}{#2}}} 
\newcommand{\gkpSII}[2]{\ensuremath{\genfrac{\lbrace}{\rbrace}{0pt}{}{#1}{#2}}} 
\newcommand{\gkpEI}[2]{\ensuremath{\genfrac{\langle}{\rangle}{0pt}{}{#1}{#2}}} 
\newcommand{\gkpEII}[2]{\ensuremath{\left\langle\genfrac{\langle}{\rangle}{
            0pt}{}{#1}{#2}\right\rangle}} 
\newcommand{\Iverson}[1]{\ensuremath{\left[#1\right]_{\delta}}}

\newcommand{\genfracSeq}[2]{\ensuremath{\genfrac{\lVert}{\rVert}{0pt}{}{#1}{#2}}} 

\newcommand{\sref}[1]{Section \ref{#1}} 
\newcommand{\chref}[1]{Chapter \ref{#1}}

\DeclareMathOperator{\EGF}{EGF}

\DeclareMathOperator{\guess}{guess} 
 
\DeclareMathOperator{\Poly}{Poly} 
\DeclareMathOperator{\RS}{RS} 
\DeclareMathOperator{\RemSeq}{RemSeq} 
\DeclareMathOperator{\Log}{Log}

\DeclareMathOperator{\Binom}{Binom}

\newcommand{\cf}[0]{\textit{c.f.}\ }

\newcommand{\quotetext}[1]{``#1''} 
\newcommand{\quoteemph}[1]{\quotetext{\emph{#1}}}

\newcommand{\OEISPlain}{\emph{Online Encyclopedia of Integer Sequences}} 
\newcommand{\OEIS}{\OEISPlain\ } 

\newcommand{\ttemph}[1]{\url{#1}} 
\renewcommand{\ttemph}[1]{\texttt{#1}} 

\newcommand{\MmPlain}{\emph{Mathematica}} 
\newcommand{\Mm}{\MmPlain\ } 
\newcommand{\FSFFnPlain}{\ttemph{FindSequenceFunction}} 
\newcommand{\FSFFn}{\FSFFnPlain\ } 
\newcommand{\GFUNPlain}{\ttemph{gfun}} 
\newcommand{\GFUN}{\GFUNPlain\ } 
\newcommand{\GPSFPkgNamePlain}{\ttemph{GuessPolySequenceFormulas.m}} 
\newcommand{\GPSFPkgName}{\GPSFPkgNamePlain\ } 
\newcommand{\GSDPkgNamePlain}{\ttemph{GuessSequenceData.m}} 
\newcommand{\GSDPkgName}{\GSDPkgNamePlain\ } 
\newcommand{\GPSFPkgGuessFnNamePlain}{\ttemph{GuessPolynomialSequence}} 
\newcommand{\GPSFPkgGuessFnName}{\GPSFPkgGuessFnNamePlain\ }

\begin{document}

%\title{Coffee Consumption of Graduate Students \\
%       Trying to Finish Dissertations}
\title{%GuessPolySequenceFormulas.m: \\ 
       A Computer Algebra Package for \\ 
       Polynomial Sequence Recognition}

\author{Maxie Dion Schmidt \\ \url{maxieds@gmail.com}}

\department{Computer Science}
\msthesisVII
\advisor{Professor Roy Campbell}
\degreeyear{2014}
\committee{Professor Roy Campbell}

\maketitle

\begin{abstract}

The software package developed in the thesis 
research implements functions for the 
intelligent guessing of polynomial sequence formulas based on 
user--defined expected sequence factors of the input coefficients. 
We present a specialized hybrid approach to 
finding exact representations for polynomial sequences 
that is motivated by the need for an automated procedures to discover the 
precise forms of these sums based on user guidance, or intuition, 
as to special sequence factors present in the formulas. 
In particular, the package combines the user input on the 
expected special sequence factors in the polynomial coefficient formulas 
with calls to the existing functions as subroutines that then 
process formulas for the remaining sequence terms 
already recognized by these packages. 

The factorization--based approach to polynomial sequence recognition 
is unique to this package and allows the search functions to find 
expressions for polynomial sums 
involving Stirling numbers and other special triangular sequences that 
are not readily handled by other software packages. 
In contrast to many other sequence recognition and summation software, 
the package not provide an 
explicit proof, or certificate, for the correctness of these 
sequence formulas -- only computationally guided educated 
guesses at a complete identity generating the sequence over all $n$. 
The thesis contains a number of concrete, working examples of the 
package that are intended to both demonstrate its usage and to 
document its current sequence recognition capabilities. 
\end{abstract}

\tableofcontents
\listoffigures
\chapter*{List of Symbols and Notation}

\section*{Notation and Conventions} 

\begin{symbollist}[0.7in]
     
     \item[$\mathbb{N}$] 
     The set of natural numbers, $\mathbb{N} = \{0,1,2,3,4,\ldots\}$. 
     \item[$\mathbb{Z}^{+}$] 
     The set of positive integers, $\mathbb{Z}^{+} = \{1,2,3,4,\ldots\}$. 
     \item[$\mathbb{Z}\lbrack x\rbrack$] 
     The ring of polynomials in $x$ with coefficients in the 
     integers, $\mathbb{Z}$. 
     \item[$\mathbb{Q}\lbrack x\rbrack$] 
     The ring of polynomials in $x$ with coefficients in the 
     rational numbers, $\mathbb{Q}$. 
     \item[$\mathbb{K}\lbrack x\rbrack$] 
     The ring of polynomials in $x$ with coefficients over the 
     field $\mathbb{K}$. 
     \item[$D_z^{(j)}$] 
     The derivative, or differential, operator with respect to $z$, i.e., where 
     $D_z^{(j)}[F(z)] \equiv F^{(j)}(z)$ 
     denotes the $j^{th}$ derivative of $F(z)$, 
     provided that the $j^{th}$ derivative of the function exists. 
     %\item[$\mathbb{C}$] 
     %The set of complex numbers of the form $a+b\imath$ for real numbers %real--valued 
     %$a,b \in \mathbb{R}$. 
     %\item[$\Re(s)$] 
     %The real part of a complex number $s \in \mathbb{C}$. 
     \item[$\binom{n}{k}$] 
     The binomial coefficients. 
     \item[$\gkpSI{n}{k}$] 
     The unsigned Stirling numbers of the first kind, 
     also denoted by $(-1)^{n-k} s(n, k)$. 
     \item[$\gkpSII{n}{k}$] 
     The Stirling numbers of the second first kind, 
     also denoted by $S(n, k)$. 
     \item[$\gkpEI{n}{k}$] 
     The first--order Eulerian numbers. 
     \item[$\gkpEII{n}{k}$] 
     The second--order Eulerian numbers. 
     \item[$H_n^{(r)}$] 
     The $r$--order harmonic numbers, $H_n^{(r)} := \sum_{k=1}^{n} k^{-r}$, %and 
     where the first--order harmonic numbers are denoted in the 
     shorthand notation, $H_n \equiv H_n^{(1)}$. 
     \item[$B_n$] 
     The Bernoulli numbers. 

\end{symbollist}

\mainmatter

\chapter{Introduction} 

\section{Background and Motivation} 
\label{Chapter_Intro_Section_BGMotivation} 

The form of composite sequences involving the 
Stirling numbers of the first and second kinds are common in many applications. 
The Stirling number triangles arise naturally in 
formulas involving sums of factorial functions and in the 
symbolic, polynomial expansions of binomial coefficients and other 
factorial function variants. 
The Stirling and Eulerian number triangles also both frequently occur in 
applications involving finite sums and generating functions 
over non--negative powers of integers. 
These applications include finding closed--form expressions and formulas for 
generating functions over polynomial multiples of an arbitrary sequence. 

\subsection{Example I: Computing Derivatives of 
               Stirling Number Generating Functions} 

If $p, k \in \mathbb{N}$, the 
following modified series for the generating functions for 
polynomial multiples of the unsigned Stirling numbers of the first kind,
denoted by $\gkpSI{n}{k}$, result in the expansions %of 
\begin{align} 
\label{eqn_ModPolyPowGFs_of_the_S1StirlingNumbers-stmts_v1} 
\sum_{n=0}^{\infty} n^{k} \cdot \gkpSI{n}{p} \frac{z^n}{n!} & = 
     \sum_{j=0}^{k} \gkpSII{k}{j} z^j \cdot D_z^{(j)}\left[ 
     \frac{(-1)^{p}}{p!} \cdot \Log(1-z)^{p} 
     \right] \\ 
\notag 
\sum_{n=0}^{\infty} n^{k} \cdot \gkpSI{n+1}{p+1} \frac{z^n}{n!} & = 
     \sum_{j=0}^{k} \gkpSII{k}{j} z^j \cdot D_z^{(j)}\left[ 
     \frac{(-1)^{p}}{p!} \cdot \frac{\Log(1-z)^{p}}{(1-z)} 
     \right], 
\end{align} 
where the derivative operator, $D_z^{(j)}$, 
denotes the $j^{th}$ derivative with respect to $z$ of its input and the 
\emph{Stirling numbers of the second kind} are denoted by $\gkpSII{n}{k}$. 
Given enough familiarity with the %triangles of the 
Stirling numbers of the first kind and some trial and error, 
formulas for each of the $j^{th}$ derivatives involved in the 
expansions of \eqref{eqn_ModPolyPowGFs_of_the_S1StirlingNumbers-stmts_v1} 
are obtained by extrapolation 
from the first several cases of $j \in \mathbb{N}$ to obtain the 
finite sums 
\begin{align} 
\label{eqn_Dzj_S1pzGF_exp_form_v1} 
D_z^{(j)}\left[\frac{(-1)^{p}}{p!} \cdot \Log(1-z)^{p}\right] & = 
     \sum_{i=0}^{j} \gkpSI{j}{i} \frac{(-1)^{p-i}}{(p-i)!} \cdot 
     \frac{\Log(1-z)^{p-i}}{(1-z)^{j}} \\ 
\notag 
D_z^{(j)}\left[\frac{(-1)^{p}}{p!} \cdot \frac{\Log(1-z)^{p}}{(1-z)} 
     \right] & = 
     \sum_{i=0}^{j} \gkpSI{j+1}{j+1-i} \frac{(-1)^{p+j-i}}{(1-z)^{j+1}} 
     \cdot \frac{\Log(1-z)^{p-j+i}}{(p-j+i)!}, 
\end{align} 
where the formulas in \eqref{eqn_Dzj_S1pzGF_exp_form_v1} 
may be regarded as polynomials in the %logarithm 
function, $\Log(1-z)$. 
A proof of the correctness of these formulas is then later obtained 
formally by induction on $j$. 

\subsection{Example II: A More Challenging Application} 

A more challenging, and less straightforward, 
example arises in attempting to find an exact, 
closed--form representation for the expansions of the 
ordinary generating function for the Stirling number sequence variant, 
$S_k^{(d)}(n)$, defined as in \eqref{eqn_Skdn_seq_sum_def_v1} 
with respect to each fixed $d, k \in \mathbb{Z}^{+}$. 
\begin{align} 
\label{eqn_Skdn_seq_sum_def_v1} 
S_k^{(d)}(n) & := \sum_{j=1}^{n} 
     \binom{n}{j} \gkpSI{j+1}{k+1} \frac{(-1)^j}{j! \cdot (j+d)} \cdot 
     \frac{(n+d)!}{n!}, 
\end{align} 
The first few examples of the ordinary generating function, 
$\widetilde{S}_k^{(d)}(z)$, over $d \geq 2$ for the sequence defined by 
\eqref{eqn_Skdn_seq_sum_def_v1} are provided for reference as follows: 
\begin{align} 
%\notag 
\label{eqn_STildekdn_OGFs_first_cases_listings-stmts_v1} 
\widetilde{S}_k^{(2)}(z) & = 
     -\frac{\Log\left(1-z\right)^{k-1}}{(1-z)^2} \left[ 
     \frac{z}{(k-1)!}+\frac{(-1+z) \Log\left(1-z\right)}{k!}
     \right] \\ 
\notag 
\widetilde{S}_k^{(3)}(z) & =  
     \frac{\Log\left(1-z\right)^{k-2}}{(1-z)^3} \left[ 
     \frac{z^2}{(k-2)!}+\frac{z (-4+3 z) 
     \Log\left(1-z\right)}{(k-1)!}+\frac{\left(2-4 z+2 z^2\right) 
     \Log\left(1-z\right)^2}{k!} 
     \right] \\ 
\notag 
\widetilde{S}_k^{(4)}(z) & =  
     -\frac{\Log\left(1-z\right)^{k-3}}{(1-z)^4} \Biggl[ 
     \frac{z^3}{(k-3)!}+\frac{z^2 (-9+6 z) 
     \Log\left(1-z\right)}{(k-2)!}+\frac{z \left(18-27 z+11 
     z^2\right) \Log\left(1-z\right)^2}{(k-1)!} \\ 
\notag 
   & \phantom{=-\frac{\Log\left(1-z\right)^{k-3}}{(1-z)^4} \Biggl[\ } + 
     \frac{\left(-6+18z-18 z^2+6 z^3\right) \Log\left(1-z\right)^3}{k!}
     \Biggr] \\ 
\notag 
\widetilde{S}_k^{(5)}(z) & =  
     \frac{\Log\left(1-z\right)^{k-4}}{(1-z)^5} \Biggl[ 
     \frac{z^4}{(k-4)!}+\frac{z^3 (-16+10 z) 
     \Log\left(1-z\right)}{(k-3)!}+\frac{z^2 \left(72-96 z+35 
     z^2\right) \Log\left(1-z\right)^2}{(k-2)!} \\ 
\notag 
   & \phantom{=\ } + 
     \frac{z \left(-96+216 z-176 z^2+50 z^3\right) 
     \Log\left(1-z\right)^3}{(k-1)!}+\frac{\left(24-96 z+144 z^2-96 
     z^3+24 z^4\right) \Log\left(1-z\right)^4}{k!}
     \Biggr]. 
\end{align} 
Based observations of the first several cases of these generating functions 
in \eqref{eqn_STildekdn_OGFs_first_cases_listings-stmts_v1}, 
we rewrite the expansions of these generating functions as the sum %in the form of 
\begin{align} 
\label{eqn_Skp1sn_OGF_gmdz_seq_sum_def_exp_stmt_v1} 
\widetilde{S}_k^{(d)}(z) & := 
     \sum_{n=0}^{\infty} S_{k}^{(d)}(n) z^n = \left( 
     \frac{(-1)^{d-1} \cdot \Log(1-z)^{k+1-d}}{(1-z)^{d}} 
     \right) \times %\cdot 
     \sum_{m=0}^{d-1} \frac{\Log(1-z)^{m} \cdot z^{d-1-m}}{(k+1-d+m)!} \cdot 
     g_m^{(d)}(z). 
\end{align} 
It is clear from examining the sequence data in 
\eqref{eqn_STildekdn_OGFs_first_cases_listings-stmts_v1} that the 
formulas for the polynomials, $g_m^{(d)}(z)$, specified in %the form of 
\eqref{eqn_Skp1sn_OGF_gmdz_seq_sum_def_exp_stmt_v1} 
involve a sum over factors of the Stirling numbers of the first kind and the 
binomial coefficients. 
However, finding the precise sequence inputs in the formula for these 
polynomials with the correct corresponding multiplier terms in the sum is 
not immediately obvious from the first few example cases in 
\eqref{eqn_STildekdn_OGFs_first_cases_listings-stmts_v1}. 
%and requires some trial and error using \MmPlain. 
%% 
We then proceed forward seeking a 
formula for the polynomials, $g_m^{(d)}(z)$, in the 
general template form of %\eqref{eqn_gmdz_example_seq_searcg_exp_v2} 
\begin{equation} 
\label{eqn_gmdz_example_seq_searcg_exp_v2} 
g_m^{(d)}(z) = \sum_{i} S_1(\cdot, \cdot) \cdot \Binom(\cdot, \cdot)%^{2} 
     \times \RemSeq_1(i) \cdot \RemSeq_2(m+m_0-i) \times z^i, 
%g_m^{(d)}(z) = \sum_{i} \gkpSI{\cdot}{\cdot} \binom{\cdot}{\cdot} 
%     \times \RemSeq_1(i) \RemSeq_2{m-i} \times z^i. 
\end{equation} 
where the functions $S_1(\cdot, \cdot)$ and $\Binom(\cdot, \cdot)$ denote the 
Stirling numbers of the first kind and binomial coefficients, respectively, 
each over some unspecified index inputs to these sequence functions. 

After a few hours of frustrating trial and error with \MmPlain, 
we finally arrive at a formula for these polynomials in the form of 
\begin{align} 
\label{eqn_gmdz_seq_sum_def_v1} 
g_m^{(d)}(z) & = \sum_{i=0}^{m} \gkpSI{d-m+i}{d-m} \binom{d-1}{m-i}^2 
     (-1)^{m-i} (m-i)! \cdot z^i. %\\ 
%\notag 
%   & = TODO(reverse index). 
\end{align} 
The motivation for constructing the package routines in the thesis is to 
automate the eventual discovery of the formula in 
\eqref{eqn_gmdz_seq_sum_def_v1} based on user input as to the 
general template to the formula for these polynomials specified as in 
\eqref{eqn_gmdz_example_seq_searcg_exp_v2}. 
The automated discovery of the first pair of less 
complicated formulas given in \eqref{eqn_Dzj_S1pzGF_exp_form_v1} 
is then also possible using the package. 

%\subsection{Overview} 
%\subsection{Motivating Examples} 

\section{High--Level Description of the Package} 

The \Mm package \GPSFPkgName developed as a part of the thesis 
research implements software routines for the 
intelligent guessing of polynomial sequence formulas based on user input 
%(or user intuition) 
on the expected form of the sequence formulas. 
These functions for sequence recognition then rely on some degree of 
user intuition to correctly find closed--form formulas 
that represent the input polynomial sequence. 
The logic used to construct these routines is based 
on factorization data for the expected sequence factors of the input 
polynomial coefficient terms suggested by the user. 

The template of the polynomial sequence formulas that the 
package aims to recognize satisfies an expansion of the 
general form outlined in \eqref{eqn_pmx_gen_poly_sequence_forms-intro_stmt_v1} 
where $j, j_0 \in \mathbb{N}$, $r \in \mathbb{Z}^{+}$, and 
$x$ is some (formal) polynomial variable that may assume different forms 
in the sequences input to the package routines. 
\begin{align} 
%\notag 
\label{eqn_pmx_gen_poly_sequence_forms-intro_stmt_v1} 
\Poly_j(x) & := \sum_{i=0}^{j+j_0} \left(
     \prod_{i=1}^{r} 
     \genfracSeq{\widetilde{u}_i(j) + u_i \cdot i}{
     \widetilde{\ell}_i(j) + \ell_i \cdot i}_i 
     \right) 
     %\times U_{\guess}(j, i) 
     \times \RS_1(i) \RS_2(j+j_0-i) \cdot x^i. 
\end{align} 
The product of (triangular) sequences in the first term of 
\eqref{eqn_pmx_gen_poly_sequence_forms-intro_stmt_v1} correspond to the 
factors of the expected user--defined sequences in the 
polynomial coefficient terms where the functions 
$\widetilde{u}_i(j), \widetilde{\ell}_i(j)$ are prescribed functions 
of the sequence index $j$ and where the $u_i, \ell_i \in \mathbb{Z}$ are 
prescribed application--dependent multiples of the polynomial 
summation index $i$. 
The functions $\RS_i(\cdot)$ in the previous equation denote the 
coefficient remainder terms in these polynomial formulas, which 
should ideally correspond to comparatively simpler sequences that are 
already easily recognized by existing packages discussed in 
\sref{Section_Intro_ComparisonsOfSoftware} of the thesis below. 
These existing packages may be called as subroutines to recognize the 
sequences corresponding to the remainder terms in the 
input polynomial coefficients after the forms of the sequence factors 
expected by the user are determined by the package routines. 

\section{Plan of Attack and Aims of the Thesis Research} 
%\subsection{Plan of Attack and Aims of the Research} 

A significant part of the work for the thesis is a 
``proof of concept'' implementation of the logic to find 
polynomial sequence formulas of the form in 
\eqref{eqn_pmx_gen_poly_sequence_forms-intro_stmt_v1} based on %the 
user input of the first several terms of the sequence. 
In this implementation, the focus of the package development is in 
constructing the logic to recognize the polynomial sequence formulas in the 
form of \eqref{eqn_pmx_gen_poly_sequence_forms-intro_stmt_v1}. 
For example, 
in the absence of obvious, or known, algorithms for the factorization of an 
integer by an arbitrary sequence, the implementation of this 
part of the algorithm employed by the package is effectively treated as 
an oracle within the working source code. 
     The construction of this type of integer factorization algorithm is 
     motivated by the need for such algorithms in a more efficient 
     implementation of this package. 
     A more complete and detailed specification of these factorization 
     routines is described in the future research topics outlined in 
     \chref{Chapter_Conclusions}. 

The plan is that later, once more of the machinery for generating the proposed 
polynomial sequence formulas is in place, optimizations to the code and the 
task of finding a more efficient 
implementation to generate the factorizations of a given integer 
over multiple sequence factors may be investigated further. 
Several examples of usage of the 
sequence recognition functions in the package, including figures of the 
\Mm output, are given in \chref{Chapter_PkgUsageChapter}. 
These examples provide both the working syntax of 
\Mm programs that employ the package routines and serve to document the 
capabilities of the package current at the time of this thesis draft. 

\section{Software for Sequence Recognition} 
\label{Section_Intro_ComparisonsOfSoftware} 

\subsection{Software Packages for Sequence Recognition} 

There are a number of notable existing software packages and online resources 
geared towards guessing formulas for integer and semi--rational 
sequences based on the forms of the first few terms of a sequence. 
Notable and well--known examples include the \GFUN package for 
\emph{Maple} \citep{GFUN-PKG-DOCS}, the 
\ttemph{Rate} packages for \MmPlain \footnote{ 
     See \url{http://www.mat.univie.ac.at/~kratt/rate/rate.html}. 
} 
\citep[Appendix A]{RATE.M-PKG-DOCS}, the more recently updated 
\ttemph{Guess} package \footnote{ 
     See %also 
     \url{http://axiom-wiki.newsynthesis.org/GuessingFormulasForSequences}. 
} 
for the \ttemph{FriCAS} fork of \emph{Axiom} 
which includes enhancements to the previous packages documented in 
\citep{AXIOM-GUESS-PKG-DOCS}, and a default, built--in function, \FSFFnPlain, 
in \MmPlain. 
There are still other software packages designed to perform related 
operations aimed at recognizing auxiliary properties 
such as identifying recurrence relations and generating functions for 
sequences freely available online\footnote{ 
     See the complete list of \emph{Algorithmic Combinatorics Software} 
     on the \emph{RISC} website at 
     \url{http://www.risc.jku.at/research/combinat/software/}. 
}. 
%The list of these related software packages includes 
%several of the built--in \url{Find*} functions in \MmPlain, functions in the 
%\emph{Maple} \GFUN package described in \citep{GFUN-PKG-DOCS}, and the 
%\ttemph{Guess} and \ttemph{Stirling} packages for \MmPlain, to name a few. 
%% 
The \OEISPlain, and its email--based \ttemph{SuperSeeker} program, 
provide lookup access to a large database of integer sequences, including the 
integer--valued entries for the numerator and denominators of 
rational sequences such as the Bernoulli numbers. 
A more historical account of the development of software for 
sequence recognition is provided in 
\citep[\S 2]{AXIOM-GUESS-PKG-DOCS}. 

\subsection{Polynomial and Summation Identities Involving the Stirling Numbers} 

Notice that in the absence of some underlying structure to a sequence 
(or satisfied by its generating function), guessing functions that attempt to 
find closed--form expressions for an arbitrary sequence by 
extrapolation from the input of its first few terms are inherently 
limited in obtaining a proof to verify the correctness of the 
formulas returned. 
The routines in many software packages and in the algorithms 
described in \citep{AEQUALSB-BOOK} are able to obtain computerized proofs 
or certificates for closed--form identities obtained 
from summations involving special functions. 
The correctness of formulas obtained by packages such as \GFUN 
follow if the generating function for a sequence is \emph{holonomic}, 
or equivalently, if the sequence, say $S_n$, itself satisfies a 
\emph{$P$--recurrence} of the form 
\begin{equation} 
%\notag 
\label{eqn_holonomic_P-Rec_form-stmt_v1} 
\widehat{p}_0(n) \cdot S_{n} + \widehat{p}_1(n) \cdot S_{n+1} + \cdots + 
     \widehat{p}_r(n) \cdot S_{n+r} = 0, 
\end{equation} 
whenever $n \geq n_0$ for some fixed $n_0$, with $r \geq 1$, and where the 
coefficients, $\widehat{P}_i(n) \in \mathbb{C}\lbrack n \rbrack$, 
are polynomials for each $0 \leq i \leq r$ \citep[\S B.4]{ACOMB-BOOK}. 

As pointed out in \citep{GKP} and in \citep{STIRLING.M-PKG-DOCS}, 
unlike a number of other special sequences of interest in applications, the 
Stirling numbers are \emph{not holonomic}, or do not satisfy a homogeneous 
recurrence relation of the form in \eqref{eqn_holonomic_P-Rec_form-stmt_v1}, 
so it is reasonable to expect that 
existing software to guess sequence formulas should be at least 
somewhat limited in recognizing the exact forms of summations 
involving factors these sequences. 
The \Mm package \ttemph{Stirling.m} by M. Kauers is still able to find 
recurrences satisfied by many polynomial--like sums involving the 
Stirling and Eulerian number triangles in cases of many known and new 
summation identities. %\citep{STIRLING.M-PKG-DOCS}. 
However, the example cited in Kauers' article about the package shows 
a seemingly simple polynomial--like summation involving the 
Stirling numbers of the second kind for which a %does not satisfy a 
recurrence relation in the form of 
\eqref{eqn_holonomic_P-Rec_form-stmt_v1} fails to exist 
\citep[See \S 4]{STIRLING.M-PKG-DOCS}. 

This behavior offers some explanation as to the deficiency of 
functions like \FSFFn in recognizing formulas for sequences 
involving factors of the Stirling and Bernoulli numbers. 
We now restrict our attention to constructing software routines that 
recognize formulas for the class of polynomial sums of the form in %outlined by 
\eqref{eqn_pmx_gen_poly_sequence_forms-intro_stmt_v1} 
based on intelligent guesses as to the coefficient forms input by the user. 
In the context, the package is intended to quickly assist the user in the 
discovery of formulas that arise in practice, 
like those motivated by the examples from 
\sref{Chapter_Intro_Section_BGMotivation}, 
which we then are able to prove correct later by separate methods. 

%\subsection{Comparisions to Existing Software} 
\subsection{Comparisons of the Packages to Existing Software Routines} 

The treatment of the user--defined expected sequence factors %by the package in 
in finding formulas for input polynomial sequences 
is to consider these expected sequence terms as 
primitives in the matching formulas returned by the package. 
This treatment of the user--defined expected sequence factors 
as primitives in the search for matching formulas 
is analogous to the handling of the closed--form functions 
returned by \FSFFn in \MmPlain, such as for 
scalar or constant values, powers of (polynomials in) a variable $n$, 
factorial and gamma functions, or powers of a fixed constant, $c^n$. 

For example, acceptable formulas returned by the package for the sequence of 
generating functions for polynomial powers of $n$ may correspond to 
either of the sums in the following equation involving the 
%coefficient factors of the 
\emph{Stirling numbers of the second kind}, 
$\gkpSII{n}{k}$, or the \emph{first--order Eulerian numbers}, $\gkpEI{n}{m}$, 
when $p \in \mathbb{N}$ and $|z| < 1$ 
\citep[\S 7.4, \S 6.2]{GKP}: 
\begin{align*} 
\sum_{n=0}^{\infty} n^{p} z^n & = 
     \sum_{k=0}^{p} \gkpSII{p}{k} \frac{k! \cdot z^k}{(1-z)^{k+1}} = 
     \frac{1}{(1-z)^{p+1}} \sum_{i=0}^{p} \gkpEI{p}{i} z^{i+1}. 
\end{align*} 
The forms of sequence factors of other standard sequences, including the 
Stirling numbers of the first and second kinds, common variants of the 
binomial coefficients, $\binom{n}{k}$ and $\binom{n+m}{m}$, the 
Eulerian number triangles, and other triangular sequences of interest in 
application--specific contexts are handled similarly as primitives in the 
desired formulas output by the package routines. 

%Notice that 
The factorization--based approach to determine factors of 
expected sequences by the user in this package 
differs from the methods employed to recognize 
sequence formulas by existing sequence recognition software. 
Since this method relies on user direction %some input user directions 
as to what terms the 
sequence formulas should contain, this approach is also useful in 
determining formulas involving factors of difficult sequence forms that 
are not easily recognized by existing software packages. 
%\TODO: Thus the package may be considered to take a 
The package then employs a hybrid %approach 
of the complementary approaches noted in %cited by %pointed out by 
\citep[\S 1]{AXIOM-GUESS-PKG-DOCS} 
to the search for polynomial sequence formulas. 
Specifically, the package routines employ existing sequence recognition 
functions as a subroutine to process the reminder terms in the 
sequence after the expected special sequence factors are 
identified in the coefficients of the input polynomials. 
%% 
%Additionally, since the polynomial sequences in the form of 
%\eqref{eqn_pmx_gen_poly_sequence_forms-intro_stmt_v1} involve 
%sums that depend on the sequence index as an upper bound for the variable 
%number of polynomial terms in the sequence, the methods in this package is 
%also able to find guesses on the matching 
%sequence formulas that other packages cannot return. 

\chapter{The Guess Polynomial Sequence Function Packages} 
%\chapter{Finding Sequence Formulas with the GuessPolySequenceFormulas.m Package} 
\label{Chapter_PkgUsageChapter} 

\section{Features of the Package} 
\label{Chapter_GPSFPkg_Section_Features} 

\subsection{Overview} 

The \GPSFPkgName package is designed to recognize formulas for 
polynomial sequences in one variable based on input user observations on 
factors of the polynomial coefficients. 
The public function \GPSFPkgGuessFnName provided by the package 
attempts to perform intelligent guessing of closed--form 
summation representations for a polynomial sequence of elements, 
$p_j(x) \in \mathbb{Z}\left[x\right]$, 
%\newnote{Check $p_j(x) \in \mathbb{Z}\left[x\right]$ notations} 
based on the user insights as to the coefficient factors in the end 
formula for the sequence and the first several polynomial terms passed as 
input to the function. 
%Several %Notice that a number of 
Several particular concrete examples of uses of the package to obtain 
formulas and other identities involving the Stirling numbers and 
binomial coefficients are contained in the discussions of  
\sref{subSection_PkgUsageChapter_Full} and 
\sref{subSection_PkgUsageChapter_MoreExamples_Full} of the 
thesis below. 

%\subsection{Specification of the Package Routines} 
\subsection{Specification of the Package Routines and 
             Polynomial Sequence Formulas} 

The primary package function \GPSFPkgGuessFnName provided to the user is 
implemented in \Mm code in such a way that it is able to handle multiple 
coefficient factors of sequences expected by the user.  
The focus of the examples provided as documentation for the 
package focus on the particular cases of \quoteemph{single--factor} and 
\quoteemph{double--factor} coefficient formulas for the input polynomials. 
In particular, the package search routines are of interest in 
obtaining sequence formulas corresponding to the following 
pair of summation formulas: 
\begin{align} 
\label{eqn_Polyjx_single_factor_seq_formula_template_v1} 
\Poly_j(x) & := \sum_{i=0}^{j+j_0} 
     \genfracSeq{\widetilde{u}_1(j) + u_1 i}{\widetilde{\ell}_1(j) + \ell_1 i}_1 
     %\times U_{\guess}(j, i) 
     \times \RS_1(i) \RS_2(j+j_0-i) \cdot x^i \\ 
\label{eqn_Polyjx_double_factor_seq_formula_template_v2} 
\Poly_j(x) & := \sum_{i=0}^{j+j_0} 
     \genfracSeq{\widetilde{u}_1(j) + u_1 i}{\widetilde{\ell}_1(j) + \ell_1 i}_1 
     \genfracSeq{\widetilde{u}_2(j) + u_2 i}{\widetilde{\ell}_2(j) + \ell_2 i}_2 
     %\times U_{\guess}(j, i) 
     \times \RS_1(i) \RS_2(j+j_0-i) \cdot x^i . 
\end{align} 
The polynomials in 
\eqref{eqn_Polyjx_single_factor_seq_formula_template_v1} and 
\eqref{eqn_Polyjx_double_factor_seq_formula_template_v2} correspond to the 
single--factor and double--factor sequence formula templates, respectively. 

In the previous equations, 
$j, j_0 \in \mathbb{N}$, $u_i, \ell_i \in \mathbb{Z}$, the 
functions $\widetilde{u}_i(j)$ and $\widetilde{\ell}_i(j)$ denote some 
prescribed application--dependent functions of the sequence index, and the 
form of the remaining sequences in the polynomial coefficient 
formulas are denoted by the functions $\RS_1(\cdot)$ and $\RS_2(\cdot)$. 
The package formula search routines only currently handle linear 
functions of the summation index inputs. 
Also notice that it is assumed that at least one of the 
$RS_i(\cdot)$ sequence functions is identically one, and 
that a formula for the remaining function is either easily obtained by an 
existing sequence recognition routine such as 
\MmPlain's \FSFFn function, or may be 
later identified with a relevant entry in the \OEIS database 
\citep{OEIS-UPDATED-URL}. 

\subsection{Special Triangular Sequence Factors 
            Supported by the Package} 

The built--in subpackage \GSDPkgName included with the current package 
source code provides an \quotetext{out of the box} 
implementation of several triangular sequences of 
interest in my research and that are important in motivating the 
development of this package. 
In the current implementation of the package, these user--specified sequences 
identified in the package routines include factors of the 
(signed and unsigned) Stirling number triangles, variations of 
triangular sequences derived from the binomial coefficients, and the 
first and second--order Eulerian number triangles defined recursively as in 
\citep[\S 6.1, 6.2]{GKP} \citep[\cf \S 26.8, 26.14]{NISTHB}. 
%These correspond to the settings of the list entries 
%passed in the \url{SequenceFactors} runtime option to the sequence IDs 
%\url{s1}, \url{S1}, \url{S2}, \url{Binom}, \url{Binom2} (squares) 
Each of these respective sequences correspond to special cases of the following 
triangular recurrence relation where 
$\alpha, \beta, \gamma, \alpha^{\prime}, \beta^{\prime}, \gamma^{\prime} \in \mathbb{Z}$ 
\citep[\S 5, \S 6.1--6.2]{GKP}: 
\begin{equation} 
\notag 
\genfracSeq{n}{k} = \left(\alpha n+\beta k+\gamma\right) \genfracSeq{n-1}{k} + 
     \left(\alpha^{\prime} n+\beta^{\prime} k+\gamma^{\prime}\right) 
     \genfracSeq{n-1}{k-1} + \Iverson{n = k = 0}. 
\end{equation} 
\Mm provides several standard, built--in functions for the 
(signed) Stirling numbers of the first and second kinds, and for the 
binomial coefficients. The related %\emph{RISC} group's 
\Mm package \url{Stirling.m} authored by Manuel Kauers 
\footnote{ 
     See also the \Mm package documentation at 
     \url{http://www.risc.jku.at/research/combinat/software/ergosum/RISC/Stirling.html}. 
} 
further extends the default functions for the Stirling numbers and 
defines additional functions that implement the Eulerian number 
triangles of both orders \citep{STIRLING.M-PKG-DOCS}. 

%\subsection{Sequence Requirements} 
%\subsection{Requirements on the Form of the Input Polynomial Sequences} 
%\subsection{Some Restrictions on the Forms of the Input Polynomial Sequence} 
\subsection{Some Restrictions on the Form of the Input Polynomials} 

The package function \GPSFPkgGuessFnName is designed to find formulas for 
polynomials, $p_j(x)$, whose coefficients are integer--valued. 
The guessing function is, however, able to find formulas for semi--rational 
polynomial sequences in $\mathbb{Q}\left[x\right]$ provided that the 
first several terms of the sequence input to the function 
\GPSFPkgGuessFnName are 
normalized by a user guess function, $U_{\guess}(j, i)$, as 
described in \sref{subSection_PkgUsageChapter_SpecUserGuessFns} 
of the thesis below. 
The difficulties in handling formulas for polynomials with rational 
coefficients arise in determining strictly integer--valued factors of 
rational--valued coefficient forms. 
These implementation issues are outlined in 
\sref{subSection_FutureFeaturesChapter_RationalSeqTransforms}. 
Several suggestions for transformations that pre--process polynomials with 
rational coefficients are also suggested in the section 
as features to be implemented in a future revision of the package. 

\section{Installation and Usage of the Package Routines} 
\label{subSection_PkgUsageChapter_Full} 

\subsection{Installation}

\subsubsection{Mathematica Package Installation} 

The package requires a working installation of \Mm and 
a copy of the two source files \\ 
\GPSFPkgName and \GSDPkgName provided on the \emph{SageMathCloud} project page 
at the URL listed in the next section.
To load the package under Linux, suppose that the package files are 
located in \url{~/guess-polys-pkg}. The package is then loaded by running 
\begin{center} 
\url{<<"~/guess-polys-pkg/GuessPolySequenceFormulas.m"} 
\end{center} 
A graphical summary of the short description and revision information for the 
package is printed when the package is successfully loaded from within a 
\Mm notebook. 

\subsubsection{Sage Package Installation} 

The \Mm package routines accompanying the original Master's thesis manuscript 
from 2014 now have a counterpart in the open--source \emph{SageMath} 
application. The \emph{Python} source code to this updated software for the 
\emph{Sage} environment is freely available for non--commercial usage online at 
\url{https://github.com/maxieds/GuessPolynomialSequences}. 
Provided that there is a correctly functioning version of 
\emph{Sage}, or a user account on the \emph{SageMathCloud} servers, 
installation of the package is as simple as copying all of the 
\emph{Python}, or \texttt{*.py}, files into the current working directory 
for \emph{Sage}.

\subsection{Typical Usage} 
\label{subSection_PkgUsageChapter_TypicalUsage_and_Examples} 

The examples given in this section illustrate both the syntax and 
utility of the sequence recognition routines provided by the 
package functions. Notice that the formulas returned by the function 
are pure functions in \Mm with three ordered parameters: 
1) The polynomial sequence index; 2) An input variable that denotes the 
summation index of the formula; and 3) A parameter that specifies the 
polynomial variable. 
The graphical printing of the formula data provided in the 
figures given in this section is disabled by 
setting the runtime option \url{PrintFormulas->False}. 
The runtime option \url{FSFFunction} is also available to replace the 
default \Mm function \FSFFn by an alternate sequence handling function to 
process the formulas for the remaining sequences in the 
polynomial coefficient terms, as well as the 
formulas for the coefficient indices in the 
polynomial index $j$ and the upper index of summation in 
\eqref{eqn_Polyjx_single_factor_seq_formula_template_v1} and 
\eqref{eqn_Polyjx_double_factor_seq_formula_template_v2}. 
The most common and useful of these option settings are documented in the 
examples below and in the sections of this chapter. 

\subsubsection{Examples: Coefficient Factors Involving the 
               Stirling Numbers of the First Kind} 

\begin{figure}[h] 
     \begin{center}
     \fbox{\includegraphics{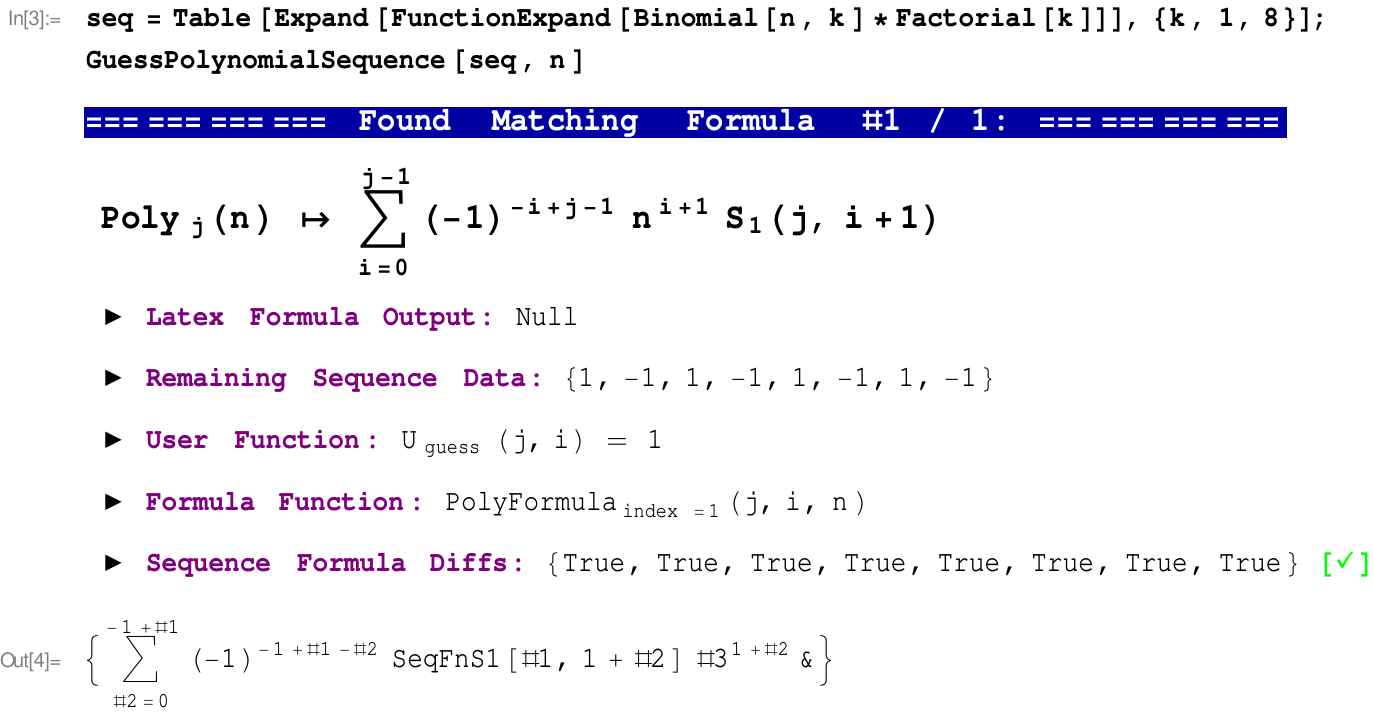}} 
          \end{center} 
     \caption{Computing a Polynomial Formula for the 
              Falling Factorial Function (\Mm)} 
     \label{figure_BinomCoeff_FallingFactFn_PolyExp_Formula_v1} 
\end{figure}

\begin{figure}[h] 
     \begin{center}
     \fbox{\includegraphics{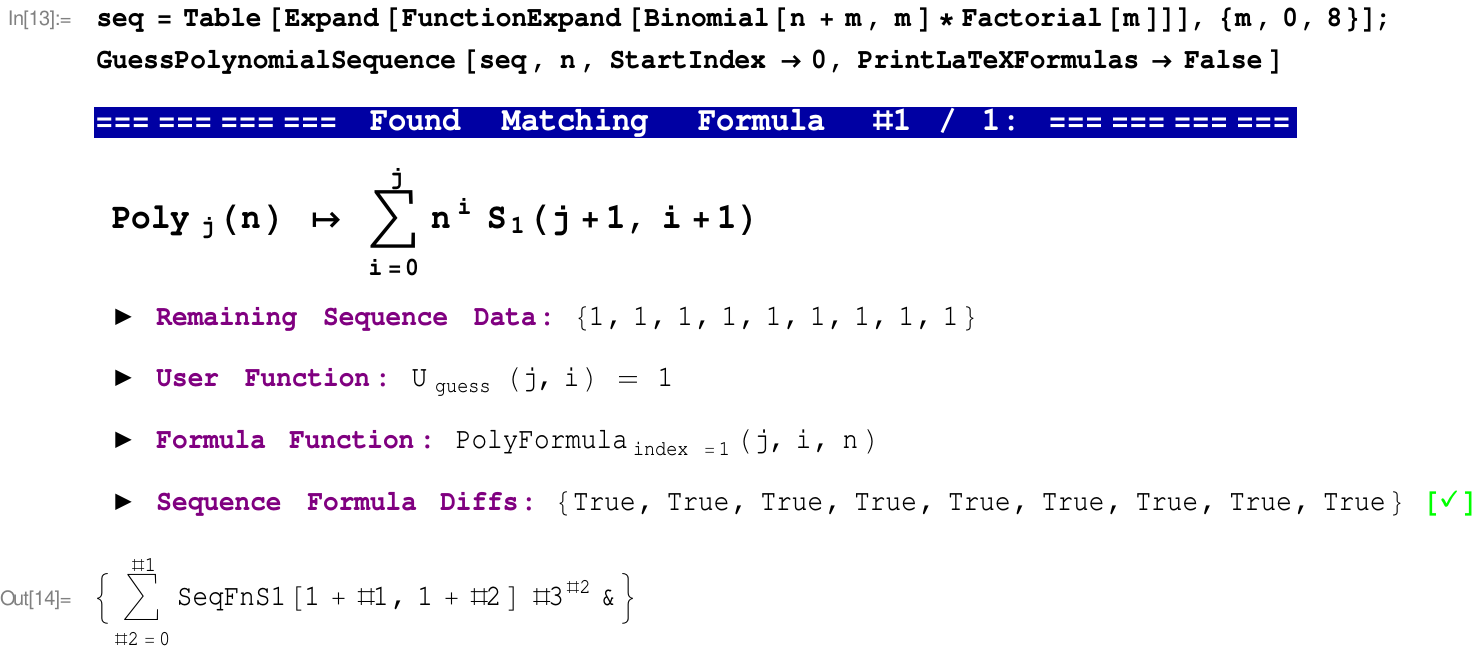}} 
     \end{center} 
     \caption{Computing a Polynomial Formula for the 
              Rising Factorial Function (\Mm)} 
     \label{figure_BinomCoeff_RisingFactFn_PolyExp_Formula_v2} 
\end{figure}

Consider the following pair of sums resulting from the 
expansions of the binomial coefficients as polynomials in $n$ 
\begin{align} 
\notag 
\binom{n}{k} & = 
     \frac{n^{\underline{k}}}{k!} = 
     \frac{1}{k!} \times n \cdot (n-1) \cdot (n-2) \cdots (n-k+1) \\ 
\label{eqn_BinomCoeff_FallingFactFn_PolyExps-stmt_v1} 
   & = 
     \frac{1}{k!} \times \sum_{i=0}^{k} \gkpSI{k}{i} (-1)^{k-i} n^{i} \\ 
\notag 
\binom{n+m}{m} & = 
     \frac{n^{\overline{m+1}}}{n \cdot m!} = 
     \frac{1}{m!} \times (n+1) \cdot (n+2) \cdots (n+m-1) \cdot (n+m) \\ 
\label{eqn_BinomCoeff_RisingFactFn_PolyExps-stmt_v2} 
   & = 
     \frac{1}{m!} \times \sum_{i=0}^{m} \gkpSI{m+1}{i+1} n^{i}, 
\end{align} 
where $n^{\underline{k}}$ denotes the \emph{falling factorial function} and 
$n^{\overline{m}}$ is the \emph{rising factorial function} %, 
in the respective expansions of the previous equations 
\citep[\S 2.6; \S 5.1]{GKP} \citep[\cf \S 26.1]{NISTHB}. 
The sums for the binomial coefficient expansions 
involving the Stirling numbers in each of 
\eqref{eqn_BinomCoeff_FallingFactFn_PolyExps-stmt_v1} and 
\eqref{eqn_BinomCoeff_RisingFactFn_PolyExps-stmt_v2} 
are known closed--form identities for the 
rising and falling factorial functions, respectively, stated in 
\citep[\S 6.1]{GKP} \citep[\cf \S 26.8]{NISTHB}. 
To see how the package can assist a user in rediscovering these 
identities, consider the respective \Mm outputs given in %the listings of 
Figure \ref{figure_BinomCoeff_FallingFactFn_PolyExp_Formula_v1} and 
Figure \ref{figure_BinomCoeff_RisingFactFn_PolyExp_Formula_v2}. 
A related example involving the polynomial sequence in 
\eqref{eqn_BinomCoeff_FallingFactFn_PolyExps-stmt_v1} is shown in 
Figure \ref{figure_BinomCoeff_FallingFactFn_PolyExp_Formula_v3}. 

\begin{figure}[h] 
     \begin{center} 
     \begin{boxedminipage}{0.87\linewidth}
     \begin{sagecommandline}
     sage: ## Falling factorial polynomials
     sage: from GuessPolynomialSequenceFunction import *
     sage: n = var('n')
     sage: poly_seq_func = lambda k: expand(simplify(binomial(n, k) * factorial(k)))
     sage: pseq_data = map(poly_seq_func, range(1, 6))
     sage: guess_polynomial_sequence(pseq_data, n, index_offset = 1);
     \end{sagecommandline}
     \end{boxedminipage}
     \end{center}

     \caption{Computing a Polynomial Formula for the 
              Falling Factorial Function (\emph{Sage})} 
     \label{figure_BinomCoeff_FallingFactFn_PolyExp_Formula_v1.2} 
\end{figure} 

\begin{figure}[h] 
     \begin{center} 
     \begin{boxedminipage}{0.87\linewidth}
     \begin{sagecommandline}
     sage: ## Rising factorial polynomials 
     sage: from GuessPolynomialSequenceFunction import *
     sage: n = var('n')
     sage: poly_seq_func = lambda m: expand(simplify(binomial(n + m, m) * factorial(m)))
     sage: pseq_data = map(poly_seq_func, range(1, 4))
     sage: guess_polynomial_sequence(pseq_data, n, index_offset = 1);
     \end{sagecommandline}
     \end{boxedminipage}
     \end{center}

     \caption{Computing a Polynomial Formula for the 
              Rising Factorial Function (\emph{Sage})} 
     \label{figure_BinomCoeff_RisingFactFn_PolyExp_Formula_v2.2} 
\end{figure} 

\subsubsection{Example: Multiple Polynomial Sequence 
               Formulas Derived from Symmetric Sequences} 

The package sequence recognition routines are able to find formulas for 
polynomials involving the \emph{binomial coefficients}, $\binom{n}{k}$, and the 
\emph{first--order Eulerian numbers}, $\gkpEI{n}{m}$. 
These sequences both have symmetry in each row of the corresponding 
triangles that satisfy the following pair of reflection identities 
where $n,k,m \in \mathbb{N}$ \citep[\S 5, \S 6.2]{GKP}: 
\begin{equation} 
\notag 
\binom{n}{k} = \binom{n}{n-k} 
     \qquad \text{ and } \qquad 
\gkpEI{n}{m} = \gkpEI{n}{n-1-m}. 
\end{equation} 
The examples given in this section demonstrate the multiple formulas 
obtained by the package for polynomial sequences involving these triangles 
that result from the coefficient symmetry noted in the forms of the 
previous equation. 

\begin{figure}[h] 
     \begin{center}
     \fbox{\includegraphics{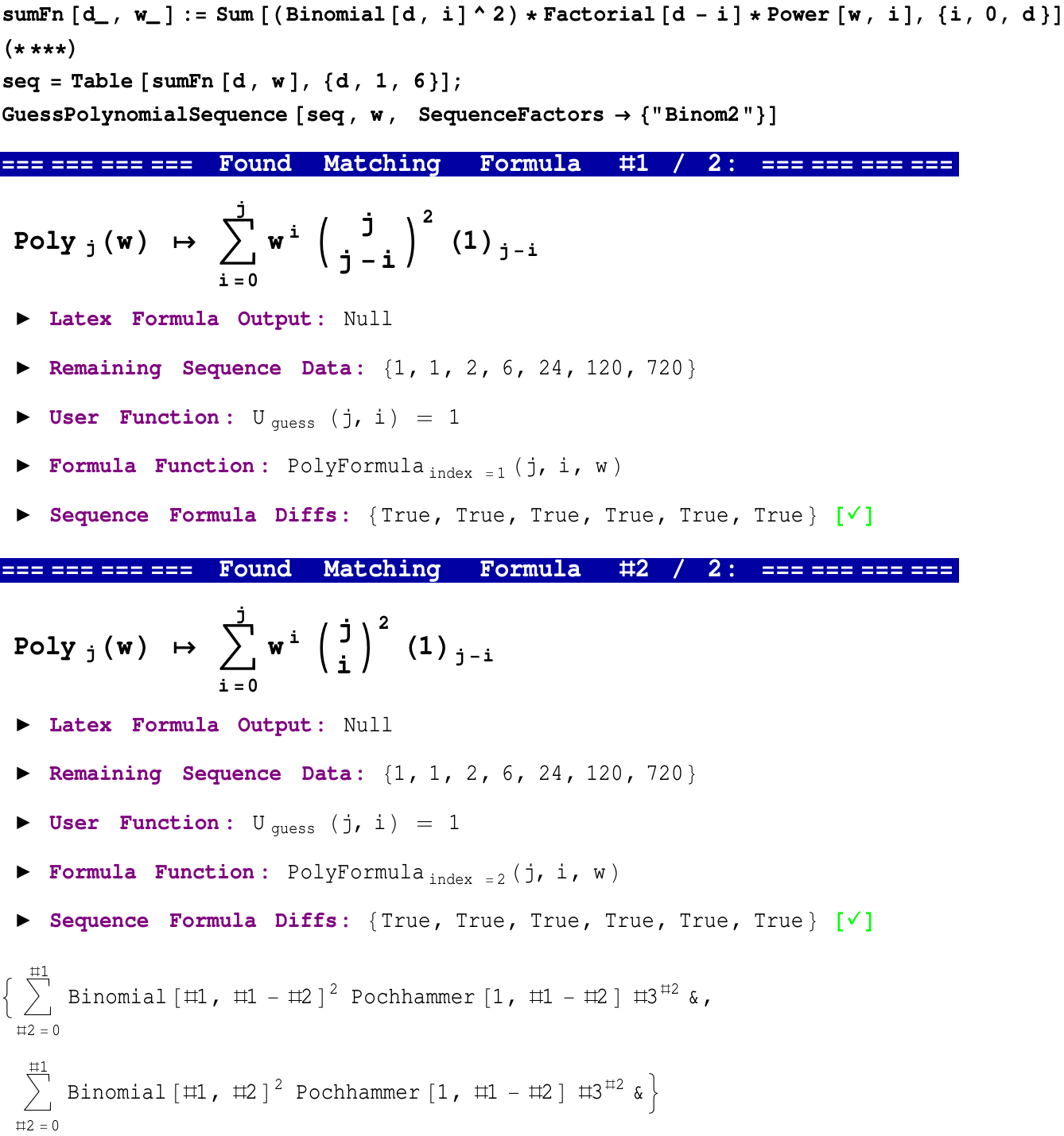}} 
     \end{center} 
     \caption{A Sum Involving Derivative Operators and 
              Squares of the Binomial Coefficients (\Mm)} 
     \label{figure_binom-derivop-formulas-v1} 
\end{figure} 

\begin{figure}[h] 
     \begin{center} 
     \begin{boxedminipage}{0.87\linewidth}
     \begin{sagecommandline}
     sage: ## Binomial squared difference operator identity
     sage: from GuessPolynomialSequenceFunction import *
     sage: n, i, w = var('n i w')
     sage: poly_seq_func = lambda d: sum((binomial(d, i) ** 2) * factorial(d - i) * (n ** i), i, 0, d)
     sage: pseq_data = map(poly_seq_func, range(1, 8))
     sage: guess_polynomial_sequence(pseq_data, n, seq_factors = ["Binom2"], index_offset = 1);
     \end{sagecommandline}
     \end{boxedminipage}
     \end{center}

     \caption{A Formula Involving Derivative Operators and 
              Squares of the Binomial Coefficients (\emph{Sage})} 
     \label{figure_binom-derivop-formulas-v1.2} 
\end{figure} 

The first example 
%shown in Figure \ref{figure_binom-derivop-formulas-v1} 
corresponds to an identity involving squares of the 
binomial coefficients and the \emph{derivative operator}, 
$D^{(j)}\left[F(z)\right] \equiv F^{(j)}(z)$, of a function, $F(z)$, 
whose $j^{th}$ derivative with respect to $z$ exists for some 
$j \in \mathbb{N}$. 
In particular, 
suppose that the function $F(z)$ denotes the 
ordinary generating function of the sequence, $\langle f_n \rangle$, and the 
function has $j^{th}$ derivatives of orders 
$j \in [0, d] \subseteq \mathbb{N}$. Then for $d \in \mathbb{Z}^{+}$, the 
generating function for the modified sequence, 
$\langle \frac{(n+d)!}{n!} f_n \rangle$, satisfies the formula 
\begin{equation} 
\label{eqn_BinomSquared_npdFactOvernFact_DerivOpIdent-stmt_v1} 
\sum_{n=0}^{\infty} \frac{(n+d)!}{n!} f_n z^n = 
     \sum_{n=0}^{\infty} (n+1) \cdots (n+d) \times f_n z^n = 
     \sum_{i=0}^{d} \binom{d}{i}^2 (d-i)! \times z^i D^{(i)}\left[F(z)\right]. 
\end{equation} 
Notice that a proof of the formula given in 
\eqref{eqn_BinomSquared_npdFactOvernFact_DerivOpIdent-stmt_v1} 
follows easily by induction on $d \geq 1$. 
A user may obtain the first several values of this sequence empirically by 
evaluating \MmPlain's \url{GeneratingFunction} for the modified 
sequence terms over the first few values of $d \geq 1$. 
Figure \ref{figure_binom-derivop-formulas-v1} shows a use of the 
package in guessing a formula for 
\eqref{eqn_BinomSquared_npdFactOvernFact_DerivOpIdent-stmt_v1} 
where the polynomial variable ($w^i$ in the figure listing) 
corresponds to the operator form of $z^i D^{(i)}$, and where the 
\emph{Pochhammer symbol}, $(1)_{j-i} \equiv (j-i)!$. 

\begin{figure}[h] 
     \begin{center}
     \fbox{\includegraphics{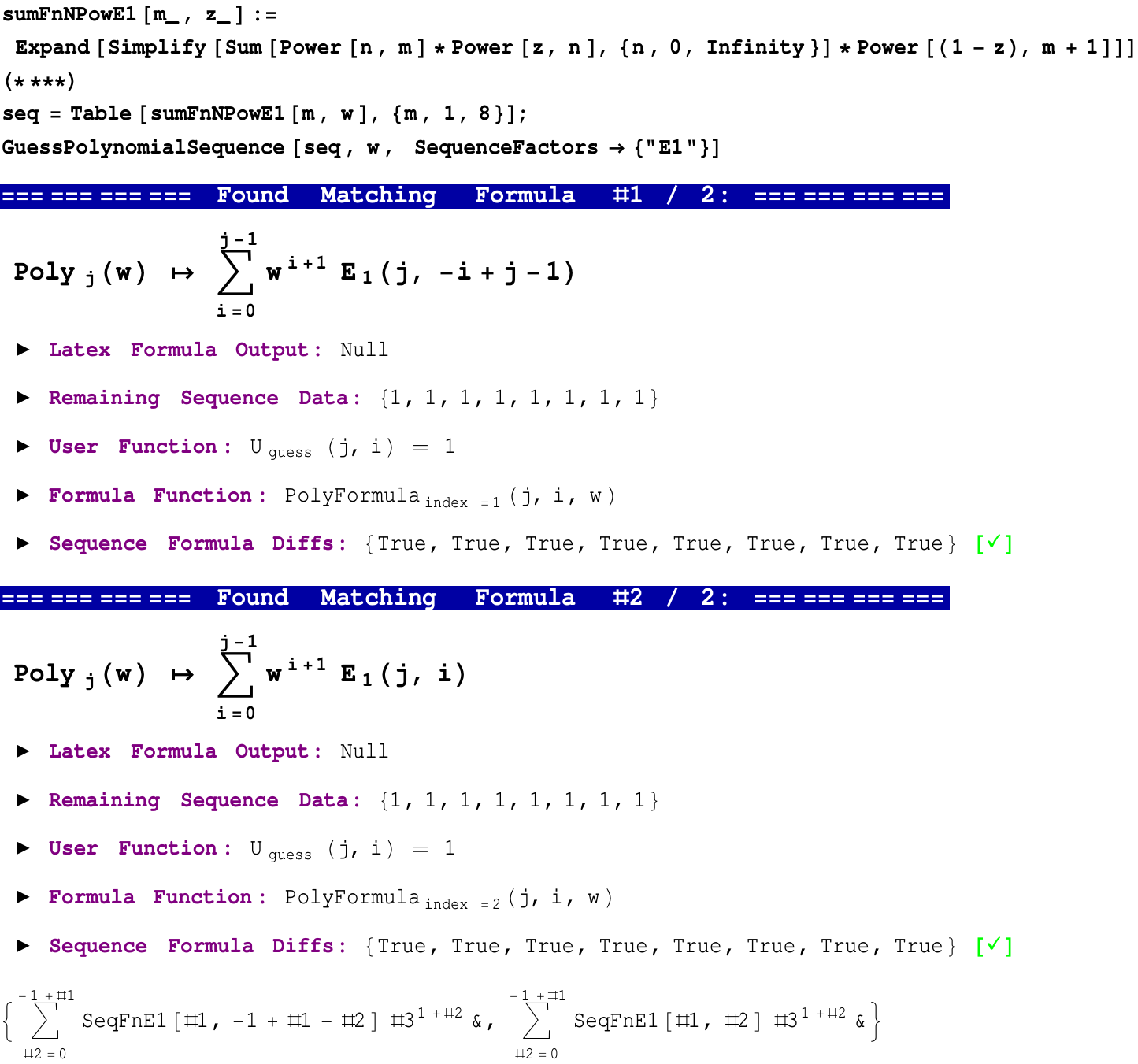}} 
     \end{center} 
     \caption{Ordinary Generating Functions of Polynomial Powers (\Mm)} 
     \label{figure_E1-geom-poly-powers-idents-v1} 
\end{figure} 

\begin{figure}[h] 
     \begin{center} 
     \begin{boxedminipage}{0.87\linewidth}
     \begin{sagecommandline}     
     sage: ## First-order Eulerian numbers in the OGFs of the polylogarithm 
     sage: ## functions, Li_{-m}(z), for natural numbers m >= 1
     sage: from GuessPolynomialSequenceFunction import *
     sage: n, z = var('n z')
     sage: poly_seq_func = lambda m: expand(factor(sum((n ** m) * (z ** n), n, 0, infinity)) * ((1 - z) ** (m + 1))).subs(z = n)
     sage: pseq_data = map(poly_seq_func, range(1, 6))
     sage: guess_polynomial_sequence(pseq_data, n, seq_factors = ["E1"], index_offset = 1);
     \end{sagecommandline}
     \end{boxedminipage}
     \end{center}

     \caption{Ordinary Generating Functions of Polynomial Powers (\emph{Sage})} 
     \label{figure_E1-geom-poly-powers-idents-v1.2} 
\end{figure} 

The next example in this section corresponds to the sequence of 
ordinary generating functions for polynomial powers of $n$, 
$\sum_{n \geq 0} n^m z^n$ for $m \in \mathbb{N}$. 
These generating functions satisfy well--known polynomial identities 
involving the \emph{Stirling numbers of the second kind}, 
$\gkpSII{n}{k}$, and the first--order Eulerian numbers, $\gkpEI{n}{m}$, 
stated as follows \cite[\cf \S 7.4]{GKP}: 
\begin{align} 
\label{eqn_PolyPowersOGFs_S2E1_sum_formulas-stmts_v1} 
\sum_{n=0}^{\infty} n^{m} z^n & = 
     \sum_{k=0}^{m} \gkpSII{m}{k} \frac{k! \cdot z^k}{(1-z)^{k+1}} = 
     \frac{1}{(1-z)^{m+1}} \sum_{i=0}^{m} \gkpEI{m}{i} z^{i+1}. 
\end{align} 
The second example cited in this section focuses on the second 
expansion in \eqref{eqn_PolyPowersOGFs_S2E1_sum_formulas-stmts_v1} 
given in terms of the Eulerian number triangle. 
Figure \ref{figure_E1-geom-poly-powers-idents-v1} 
shows the output of the package on the second polynomial sequence 
scaled by a multiple of $(1-z)^{m+1}$. 
As with the first example, the row--wise symmetry in the Eulerian number 
triangle results in the two separate formulas in the figure. 

\subsubsection{Examples: Two--Factor Polynomial Sequence Formulas} 

\begin{figure}[h] 
     \begin{center}
     \fbox{\includegraphics{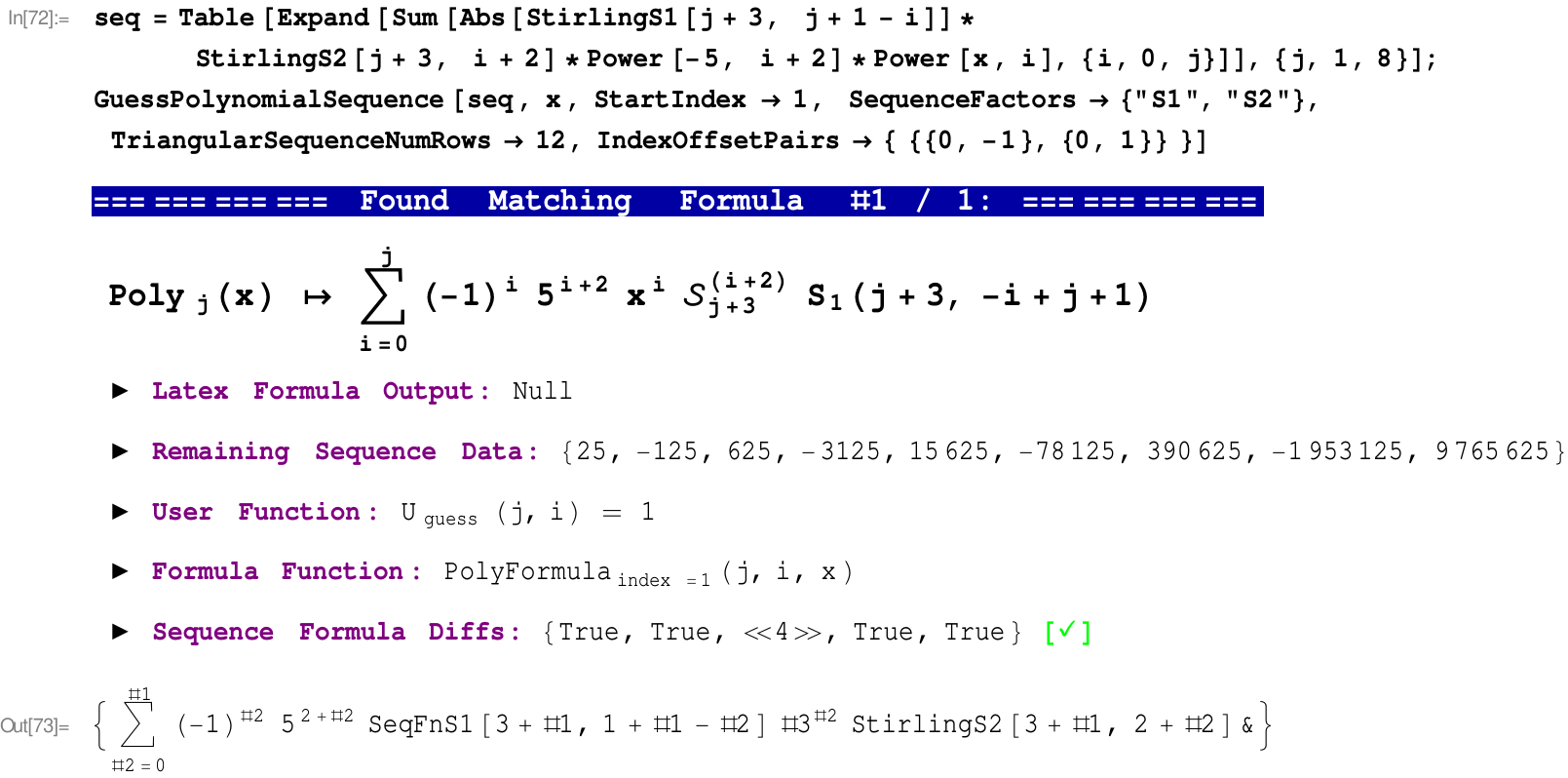}} 
     \end{center} 
     \caption{A Double--Factor Sequence Example Involving the 
              Stirling Number Triangles} 
     \label{figure_S1S2-double-factor-formula-example-v1} 
\end{figure} 

\begin{figure}[h!] 
     \begin{center}
     \fbox{\includegraphics{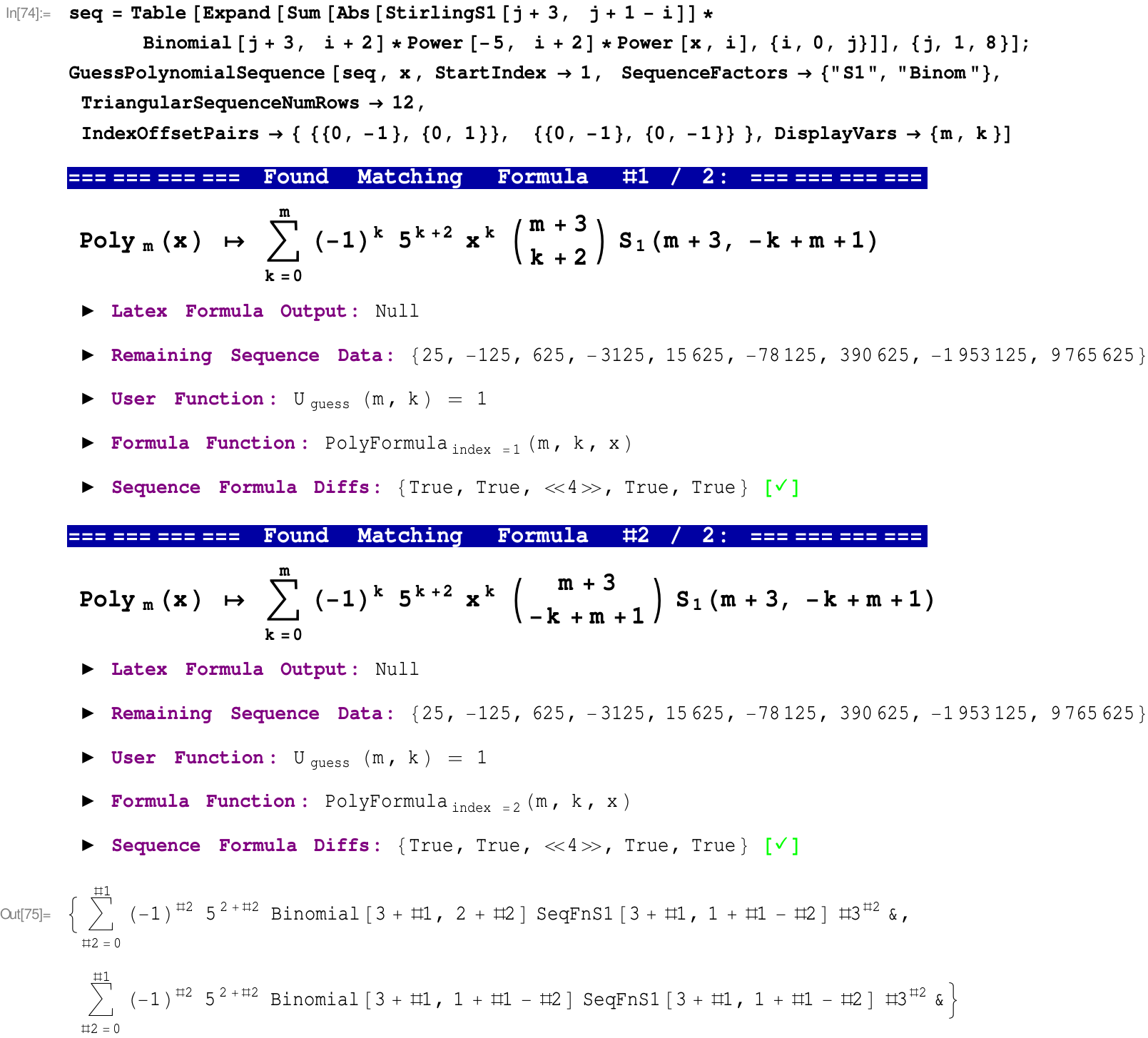}} 
     \end{center} 
     \caption{A Double--Factor Sequence Example Involving the 
              Stirling Numbers of the First Kind and the 
              Binomial Coefficients} 
     \label{figure_S1Binom-double-factor-formula-example-v2} 
\end{figure} 

The examples cited in this section correspond to the two--factor 
polynomial sequence formulas in the form of 
\eqref{eqn_Polyjx_double_factor_seq_formula_template_v2}. 
The first example given in 
Figure \ref{figure_S1S2-double-factor-formula-example-v1} 
shows the output of the package function \\ 
\GPSFPkgGuessFnName for a formula involving the 
Stirling numbers of the first and second kinds. The second example given in 
Figure \ref{figure_S1Binom-double-factor-formula-example-v2} 
shows the pair of formulas output for a sequence formula involving the 
Stirling numbers of the first kind and the binomial coefficients. 
The use of the runtime option \url{IndexOffsetPairs} in both of these 
examples is explained in more detail by 
\sref{subSection_PkgUsageChapter_Troubleshooting_PotentialIssues}. 

\subsection{User Guess Functions} 
\label{subSection_PkgUsageChapter_SpecUserGuessFns} 

The guessing routines implemented in the package rely on some intuition on the 
part of the user to determine a general template for the end formulas for 
an input polynomial sequence with coefficients over the integers. 
The user may specify an additional ``\emph{user guess function}'' that is 
employed by the package to pre--process the coefficients of the 
polynomial sequence terms passed to the function \GPSFPkgGuessFnName. 
This construction allows semi--rational, and even non--polynomial functions 
in the input variable to be processed by the package functions. 

\subsubsection{Example: A Second Formula for the Falling Factorial Function} 

\begin{figure}[h] 
     \begin{center}
     \fbox{\includegraphics{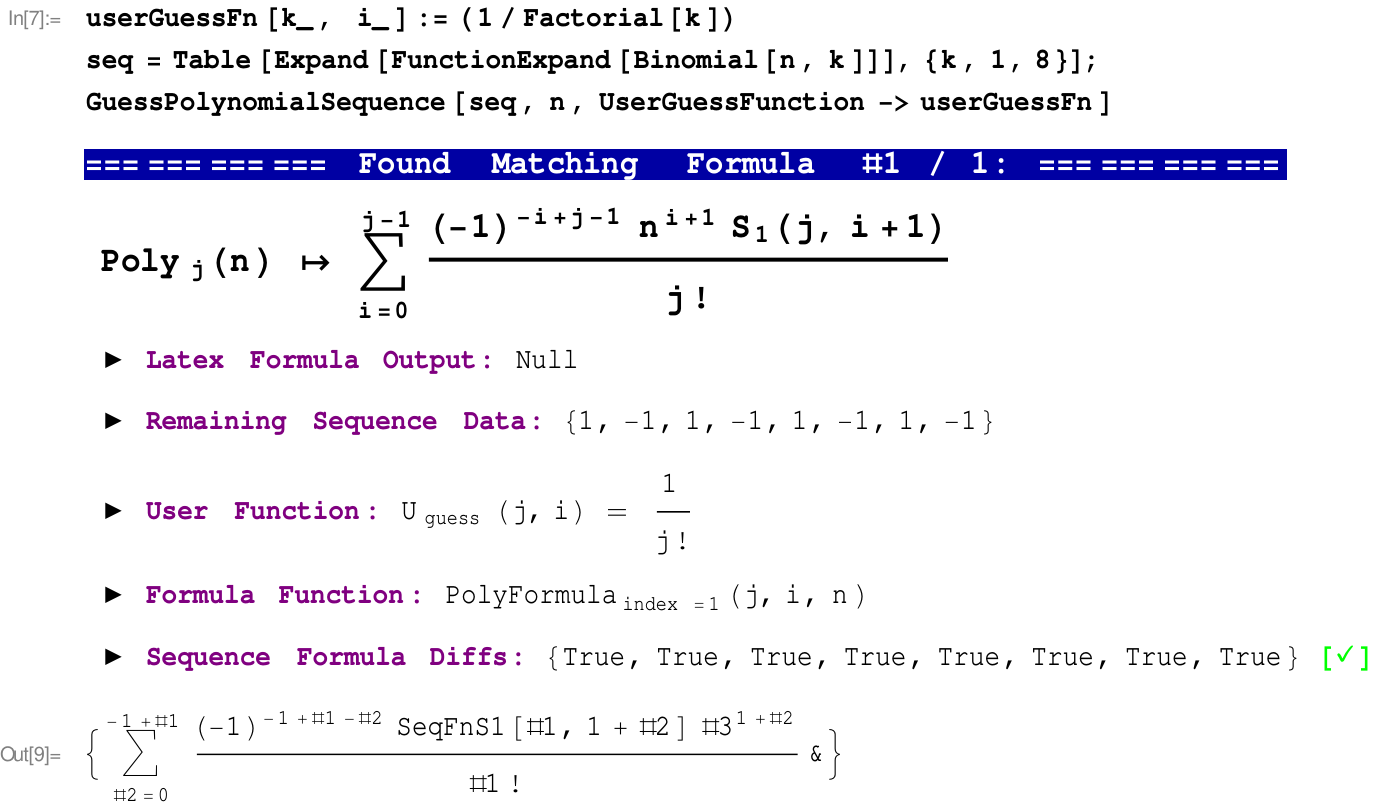}} 
     \end{center} 
     \caption{A Formula for the Falling Factorial Function 
              by a User Guess Function (\Mm)} 
     \label{figure_BinomCoeff_FallingFactFn_PolyExp_Formula_v3} 
\end{figure} 

\begin{figure}[h] 
     \begin{center}
     \begin{boxedminipage}{0.87\linewidth}
     \begin{sagecommandline}     
     sage: ## Another exponential falling factorial polynomial example 
     sage: ## with a user guess function 
     sage: from GuessPolynomialSequenceFunction import * 
     sage: n, z = var('n z')
     sage: poly_seq_func = lambda k: binomial(n, k)
     sage: user_guess_func = lambda n, k: 1 / factorial(k)
     sage: pseq_data = map(poly_seq_func, range(1, 6))
     sage: guess_polynomial_sequence(pseq_data, n, user_guess_func = user_guess_func, index_offset = 1);
     \end{sagecommandline}
     \end{boxedminipage}
     \end{center} 
     \caption{Computing a Formula for the Falling Factorial Function 
              by a User Guess Function (\emph{Sage})} 
     \label{figure_BinomCoeff_FallingFactFn_PolyExp_Formula_v3.2} 
\end{figure} 

A first example of the syntax for guessing the polynomial expansions of the 
binomial coefficient identity from 
\eqref{eqn_BinomCoeff_FallingFactFn_PolyExps-stmt_v1} is provided in 
Figure \ref{figure_BinomCoeff_FallingFactFn_PolyExp_Formula_v3}. 
Notice that this example is similar to the first form of the 
sequence formula computed by the package in 
Figure \ref{figure_BinomCoeff_FallingFactFn_PolyExp_Formula_v1}, 
except that in this case the input sequence is not normalized by a factor 
of $k$ to make the polynomial coefficients strictly integer--valued. 
A similar computation is employed to discover an analogous sum for the 
non--normalized sequence formula corresponding to the 
rising factorial function from 
Figure \ref{figure_BinomCoeff_RisingFactFn_PolyExp_Formula_v2}. 

\subsubsection{Example: An Exponential Generating Function for the 
                Binomial Coefficients} 

\begin{figure}[h!] 
     \begin{center}
     \fbox{\includegraphics[width=0.95\textwidth]{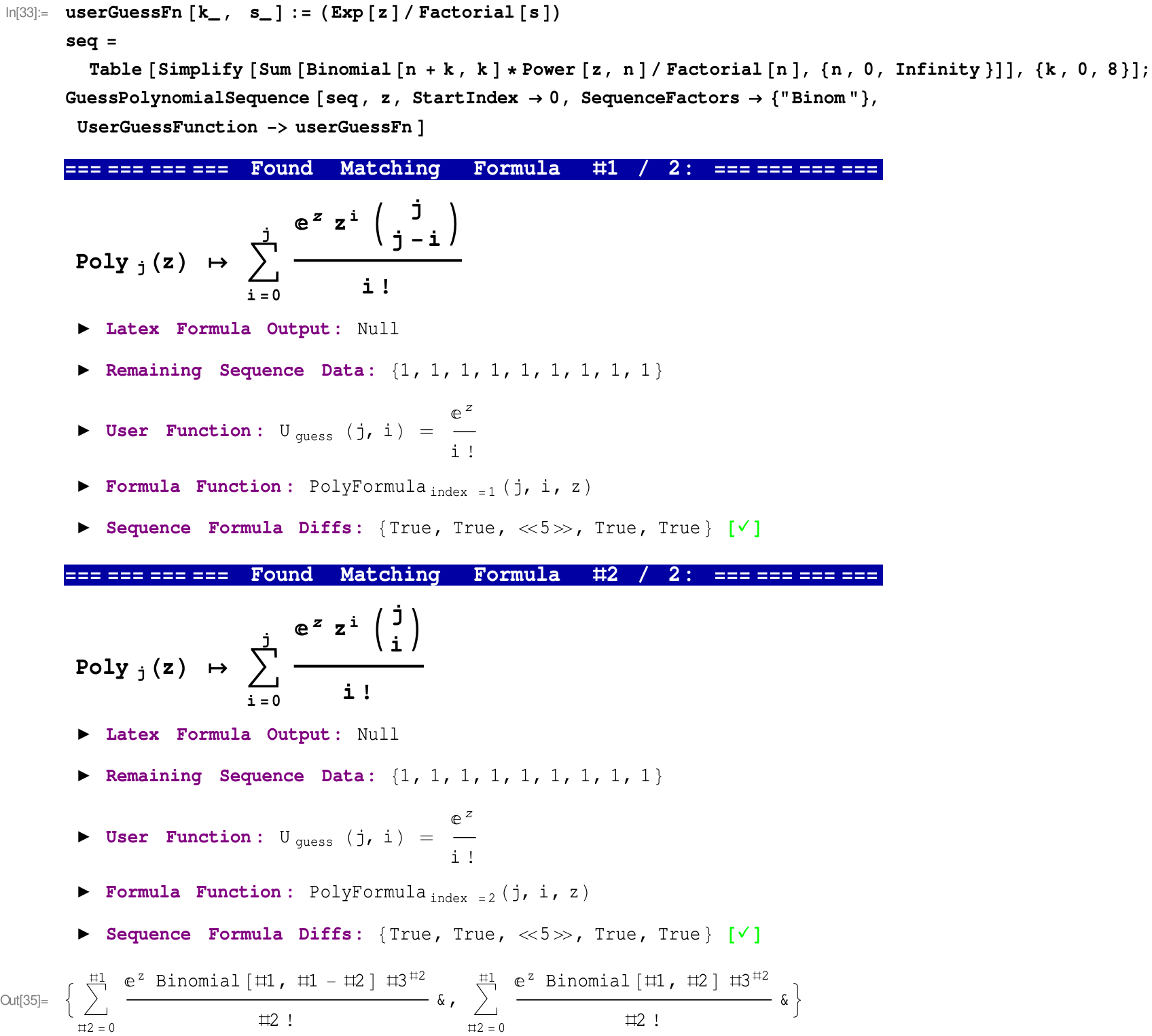}} 
     \end{center} 
     \caption{An Exponential Generating Function for the 
              Binomial Coefficients (\Mm)} 
     \label{figure_binom_EGF_formula_v1} 
\end{figure} 

\begin{figure}[h] 
     \begin{center}
     \begin{boxedminipage}{0.87\linewidth}
     \begin{sagecommandline}     
     sage: ## Exponential generating functions for the symmetric--indexed 
     sage: ## binomial coefficients
     sage: from GuessPolynomialSequenceFunction import *
     sage: n, z = var('n z')
     sage: poly_seq_func = lambda k: sum(binomial(n + k, k) * (z ** n) / factorial(n), n, 0, infinity) * exp(-z)
     sage: user_guess_func = lambda n, k: 1 / factorial(k)
     sage: pseq_data = map(poly_seq_func, range(1, 6))
     sage: guess_polynomial_sequence(pseq_data, z, seq_factors = ["Binom2"], user_guess_func = user_guess_func, index_offset = 1);
     \end{sagecommandline}
     \end{boxedminipage}
     \end{center} 
     \caption{An Exponential Generating Function for the 
              Binomial Coefficients (\emph{Sage})} 
     \label{figure_binom_EGF_formula_v1.2} 
\end{figure} 

A sequence of exponential generating functions for the symmetric form of the 
binomial coefficients, $\binom{n+m}{m}$, taken over $m \in \mathbb{N}$ 
satisfies the formula given in the following equation 
\citep[\cf \S 7.2]{GKP}: 
\begin{equation} 
%\notag 
\label{eqn_binom_EGF_formula_v1} 
%\widehat{B}_{m+1}(z) := 
\EGF_z\left(\frac{1}{(1-z)^{m+1}}\right) \equiv 
     \sum_{n=0}^{\infty} \binom{n+m}{n} \frac{z^n}{n!} \equiv 
     \sum_{s=0}^{m} \binom{m}{s} \frac{e^{z} \cdot z^s}{s!}. 
\end{equation} 
A proof of this identity is given using \emph{Vandermonde's convolution} 
identity for the binomial coefficients 
\citep[\S Table 174; \S 5.2; \cf eq. (5.22)]{GKP}. 
Figure \ref{figure_binom_EGF_formula_v1} 
shows a use of the package to guess the formula in 
\eqref{eqn_binom_EGF_formula_v1} by providing 
a user guess function that effectively removes the 
factor of $e^{z}$ in the expected formula, and that cancels out the 
coefficient factors of $1/s!$ to produce an input sequence with 
integer coefficients. 

\subsection{Troubleshooting Possible Issues} 
\label{subSection_PkgUsageChapter_Troubleshooting_PotentialIssues} 

\subsubsection{Inputting an Insufficient Number of Sequence Elements} 

There are a couple of issues that can arise in running the package 
routines when too few values of the sequence are passed to the 
\GPSFPkgGuessFnName function. The first of these is that 
\FSFFn may require a lower bound on the number of sequence values 
necessary to compute formulas for the remaining sequence terms. 
This can occur, for example, when the remaining sequence is a 
polynomial in the summation index. 
\begin{figure}[h!] 
     \begin{center}
     \fbox{\includegraphics[width=0.95\textwidth]{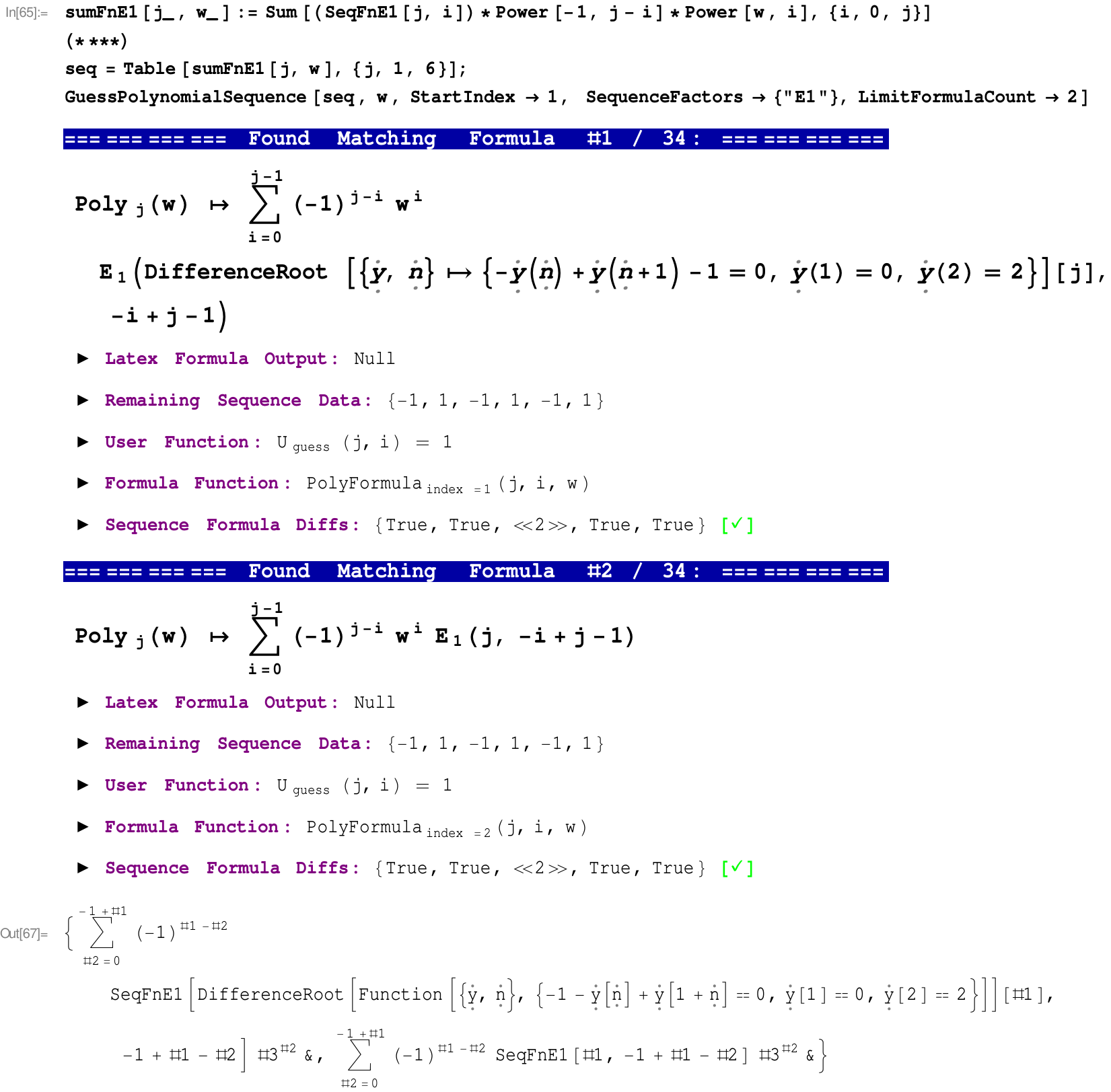}} 
     \end{center} 
     \caption{Troubleshooting an Insufficient Number of Sequence Elements} 
     \label{figure_insufficient_elements_example_v1} 
\end{figure} 
Another quirk of \MmPlain's built--in \FSFFn is that it may return a 
sequence formula matching a recurrence relation that is 
actually accurate for the few sequence elements input to the function. 
An example of this behavior is illustrated by the output given in 
Figure \ref{figure_insufficient_elements_example_v1}. 
In most cases, the problem is resolved by simply passing more 
polynomials from the sequence, usually at least $6$, but possibly 
$8$ or more elements from the sequence. 
The package is configured to warn users when less than $6$ initial terms 
are input to the function with no matching formulas. 

\begin{figure}[h!] 
     \begin{center}
     \fbox{\includegraphics[width=0.95\textwidth]{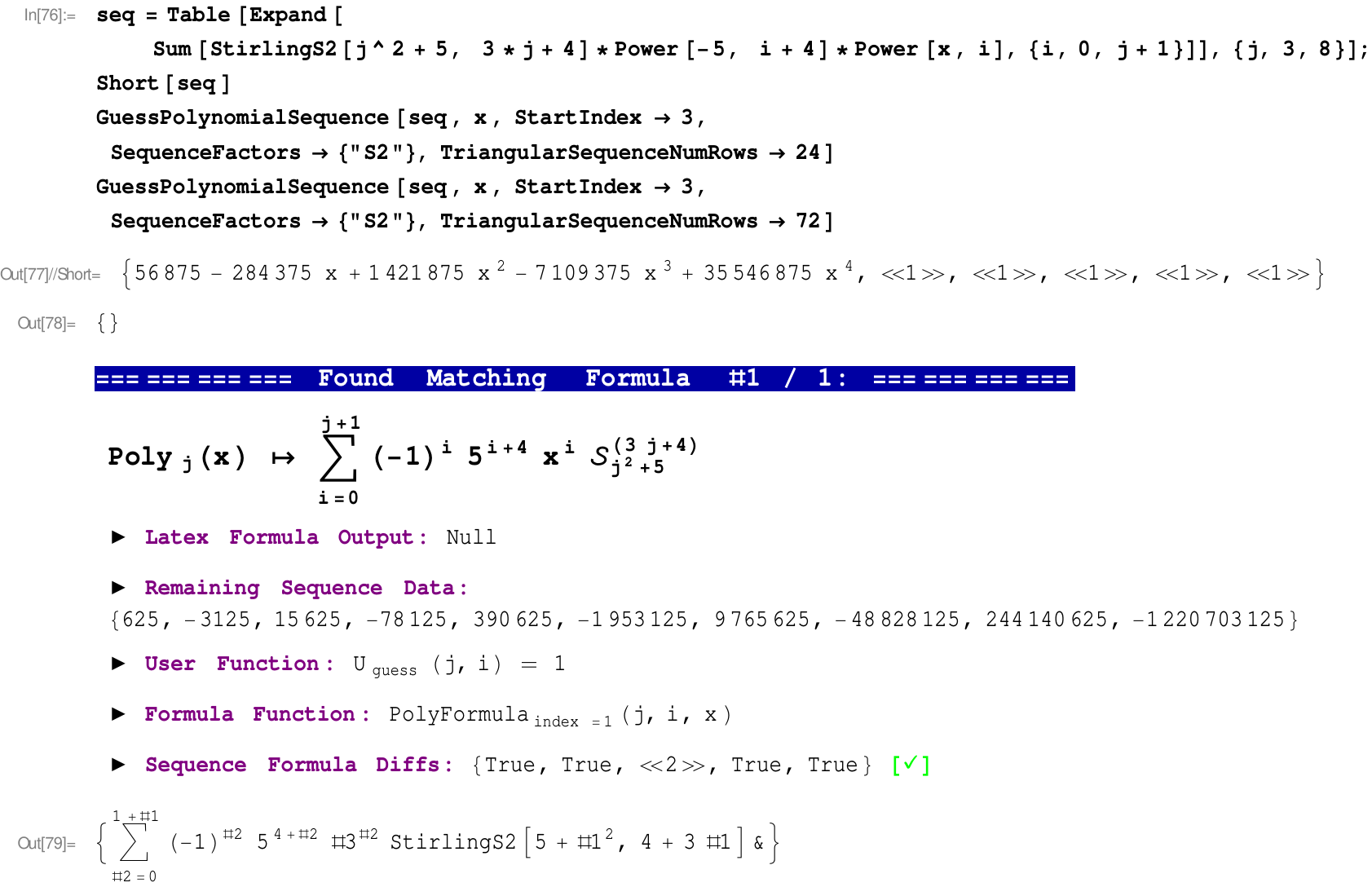}} 
     \end{center} 
     \caption{Troubleshooting Runtime Settings of the 
              \texttt{TriangularSequenceNumRows} Option} 
     \label{figure_triangle-num-rows-example-v1} 
\end{figure} 

\subsubsection{Number of Rows for the Expected Triangular Sequence Factors} 

In some cases, the package functions may not be able to obtain a formula 
for an input sequence due to an insufficient setting for the number of 
rows to consider for the expected triangular sequence factors. 
The runtime option to change the number of rows used to detect the 
factors of the expected triangular sequence is 
\url{TriangularSequenceNumRows} 
(the current default setting is \url{TriangularSequenceNumRows->12}). 
Figure \ref{figure_triangle-num-rows-example-v1} 
provides an example involving the Stirling numbers of the second kind 
where the upper index of the sequence depends quadratically on the 
polynomial index $j$. 
In this example, the package routines are unable to obtain a formula when the 
runtime option is reset to \url{TriangularSequenceNumRows->24}, but 
correctly finds the sequence formula by setting the option to the 
higher value of \url{TriangularSequenceNumRows->72}. 
Notice that choosing a significantly higher default setting for this 
option may result in much slower running times, especially if the 
expected triangular sequence factors contain a large number of $1$--valued 
entries, for example, as in the Stirling numbers of the first kind, 
binomial coefficient, and first--order Eulerian number triangles. 

\subsubsection{Handling Long Running Times with Multiple Sequence Factors} 

The package function \GPSFPkgGuessFnName is able to return 
sequence formulas in the single--factor form given in 
\eqref{eqn_Polyjx_single_factor_seq_formula_template_v1} 
in a reasonable amount of running time. 
As suggested in the double--factor sequence examples of the form in 
\eqref{eqn_Polyjx_double_factor_seq_formula_template_v2} from 
Figure \ref{figure_S1S2-double-factor-formula-example-v1} and 
Figure \ref{figure_S1Binom-double-factor-formula-example-v2}, the 
runtime option \url{IndexOffsetPairs} is needed to speed--up the 
running time for the computations involved in these sequence cases. 
The \url{IndexOffsetPairs} option is defined as a list of lists of the form 
\begin{equation} 
%\notag 
\label{eqn_IndexOffsetPairs_option_element_forms_v1} 
\mathtt{ 
     \{ \{u_1, \ell_1\}, \{u_2, \ell_2\}, \ldots, \{u_r, \ell_r\} \} 
}, 
\end{equation} 
where $r \geq 1$ denotes the expected number of sequence factors 
involved in the search for the sequence formulas by \GPSFPkgGuessFnNamePlain. 
In the examples cited in 
Figure \ref{figure_S1S2-double-factor-formula-example-v1}, 
Figure \ref{figure_S1Binom-double-factor-formula-example-v2}, and in the 
template form of \eqref{eqn_Polyjx_double_factor_seq_formula_template_v2}, the 
value of $r$ corresponds to $r := 2$. 
For a fixed choice of $r \geq 1$, 
each element of the list defined by \url{IndexOffsetPairs} passed in the 
form of \eqref{eqn_IndexOffsetPairs_option_element_forms_v1} 
corresponds to a search for a sequence formula of the form 
\begin{align} 
\notag 
%\label{eqn_Polyjx_single_factor_seq_formula_template_v1} 
\Poly_j(x) & := \sum_{i=0}^{j+j_0} \left(
     \prod_{i=1}^{r} 
     \genfracSeq{\widetilde{u}_i(j) + u_i i}{\widetilde{\ell}_i(j) + \ell_i i}_i 
     \right) 
     %\times U_{\guess}(j, i) 
     \times \RS_1(i) \RS_2(j+j_0-i) \cdot x^i. 
\end{align} 
Thus resetting the value of this option at runtime can speed--up the 
search for matching formulas in the cases of multiple expected 
sequence factors, especially 
compared to the number of index offset pairs resulting 
from the default enumeration of these pair values. 

%\section{Other Examples of Sequences Recognized by the Package} 
%\section{Other Sequence Types Recognized by the Package and Examples} 
%\section{Sequence Types Recognized by the Package} 
\section{More Examples of Polynomial Sequence Types Recognized by the Packages} 
\label{subSection_PkgUsageChapter_MoreExamples_Full} 

The examples cited in this section are intended to document further 
forms of the polynomial sequence types that the package is able to 
recognize. These examples include handling polynomial sequence formulas 
that depend on arithmetic progressions of indices, 
coefficients that contain symbolic data, and examples of sequence formulas 
obtained by the package routines when the expected sequence factors do not 
depend on the summation index, i.e., when the factors only depend on the 
polynomial sequence index. 

\begin{figure}[h!] 
     \begin{center}
     \fbox{\includegraphics[width=0.95\textwidth]{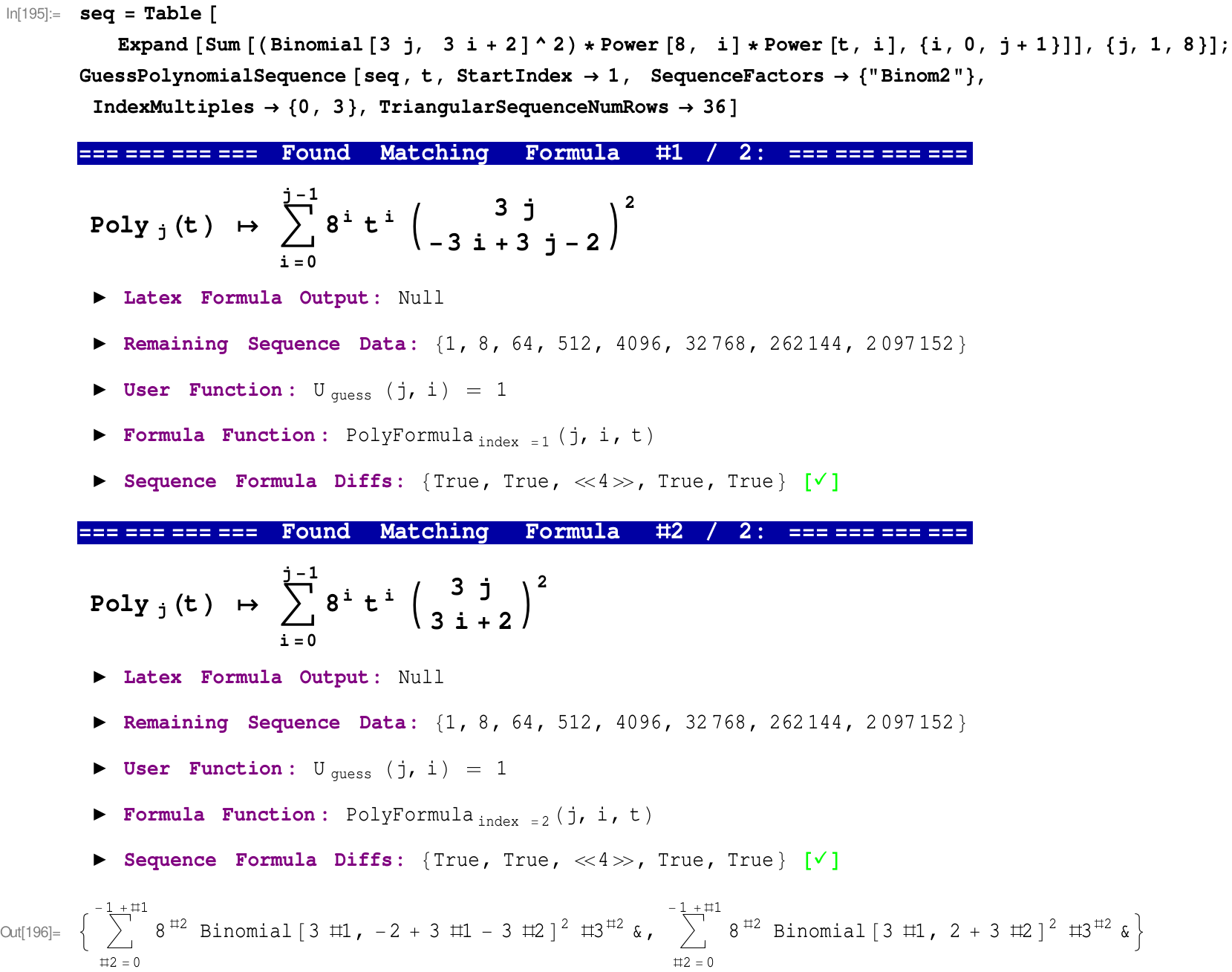}} 
     \end{center} 
     \caption{Handling Arithmetic Progressions of Indices} 
     \label{figure_arith-progs-of-coeffs-example-v1} 
\end{figure} 

\subsection{Example: Arithmetic Progressions of Coefficient Indices} 

The package function \GPSFPkgGuessFnName can be configured to search for 
sequence formulas involving arithmetic progressions of the summation 
index, $f(j)+ai$, for values besides $a := \pm 1$ by resetting the 
runtime option \url{IndexMultiples}. The default setting of this option is 
\url{IndexMultiples->{0,1}}. 
Figure \ref{figure_arith-progs-of-coeffs-example-v1} 
provides an example of recognizing sequence formulas involving squares of the 
binomial coefficients where the upper index of the triangle does not 
depend on the summation index (a setting of $a := 0$) and where the 
lower triangle index involves an arithmetic progression of the summation 
index with $a := \pm 3$. 
Related sequence formulas are recognized by setting the runtime value of this 
option to a list of test values that is some subset of the natural numbers. 
Notice that if the list of values for the option \url{IndexMultiples} 
does not contain $0$, the package routines will not find formulas like those 
given in Figure \ref{figure_arith-progs-of-coeffs-example-v1} 
where the upper index of the expected triangle factors only depends on the 
polynomial sequence index ($j$ in the figure examples). 

%\subsection{Handling of Non--Polynomial Sequences in One Variable} 
%\subsection{Handling Non--Polynomial Sequences in One Variable} 

%\subsection{Insertion of a Dummy Polynomial Variable} 
%In some situations, the user may suspect ..., even where an explicit 
%symbolic polynomial variable is missing in the formula expected by the user. 
%In the following several examples given in this section, ... 
%(see notes) 
%This procedure may facilitate the discovery of additional, non--polynomial 
%sequence formulas. 
%\subsubsection{Example: A Summation Formula for the First--Order Eulerian Numbers} 
%\subsubsection{Example: The Second--Order Eulerian Numbers} 
%\subsubsection{Example: } 
%\subsubsection{Example: } 

\begin{figure}[h!] 
     \begin{center}
     \fbox{\includegraphics{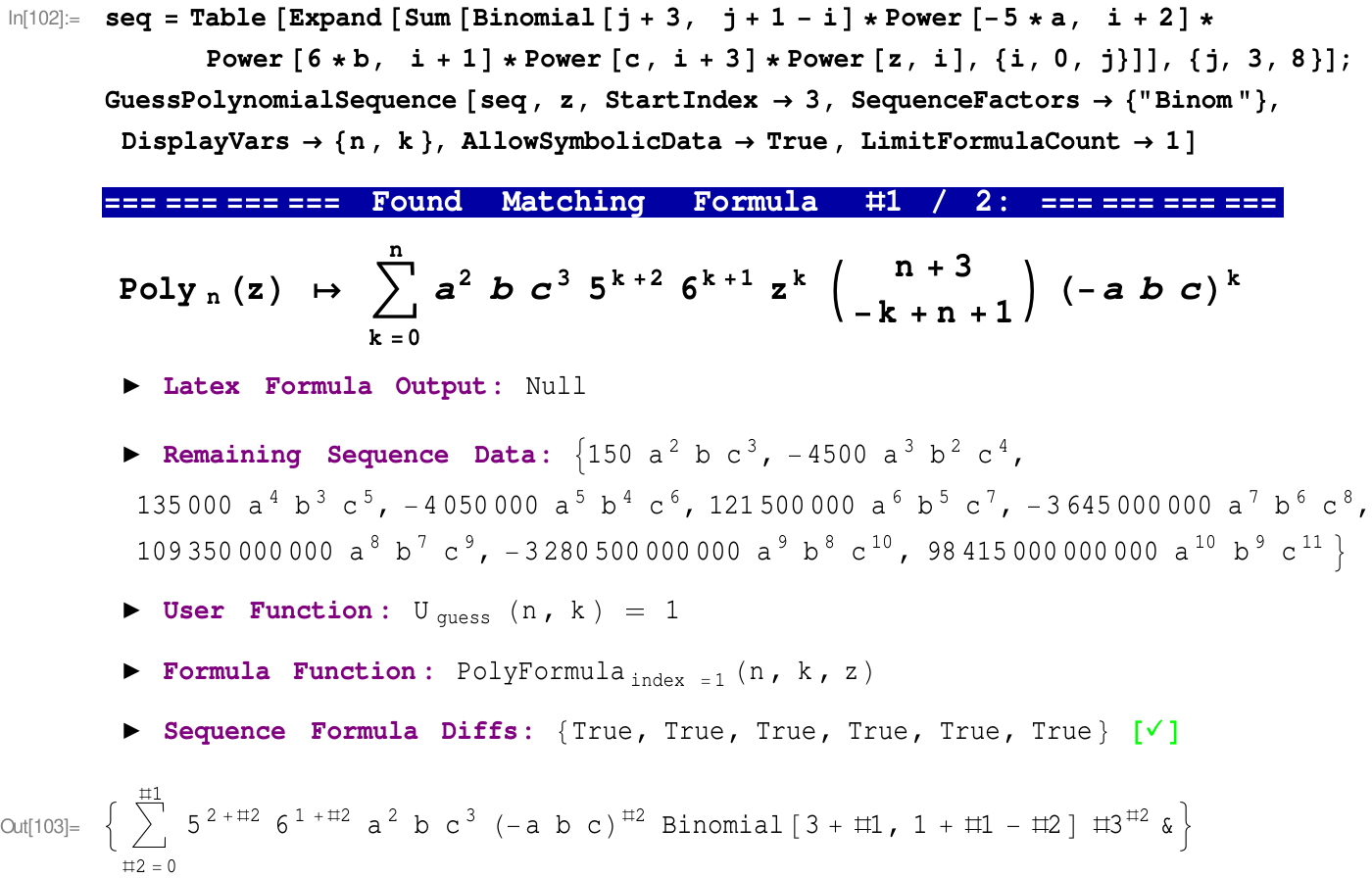}} 
     \end{center} 
     \caption{Recognizing Sequence Formulas Involving Symbolic Coefficients} 
     \label{figure_symbolic-coeff-data-example-v2} 
\end{figure} 

\begin{figure}[h!] 
     \begin{center}
     \fbox{\includegraphics{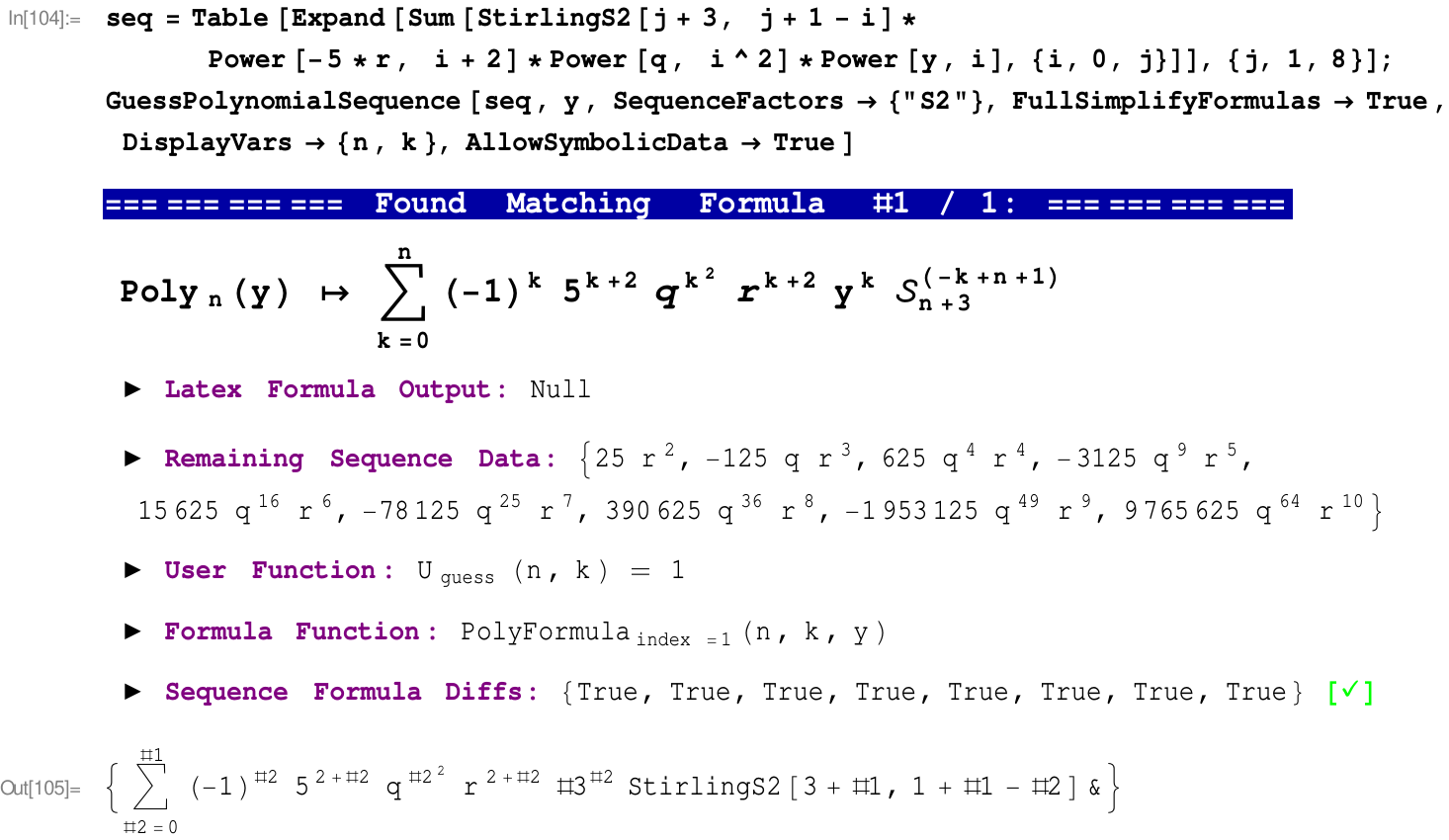}} 
     \end{center} 
     \caption{A Second Formula Involving Square Index Powers of 
              Symbolic Coefficients} 
     \label{figure_symbolic-coeff-data-example-v3} 
\end{figure} 

\subsection{Examples: Formulas Involving Symbolic Coefficient Data} 

The function \GPSFPkgGuessFnName can be configured to search for 
formulas where the input coefficients of the polynomial sequence contain 
non--numeric factors of symbolic data through the runtime option 
\url{AllowSymbolicData}. 
Figure \ref{figure_symbolic-coeff-data-example-v2} and 
Figure \ref{figure_symbolic-coeff-data-example-v3} 
provide examples of sequence formulas involving non--numeric, symbolic terms, 
named $a,b,c,q,r$, that are recognized by the package by passing \\ 
\url{AllowSymbolicData->True} to the \GPSFPkgGuessFnName function at runtime. 

\begin{figure}[h!] 
     \begin{center}
     \fbox{\includegraphics[width=0.95\textwidth]{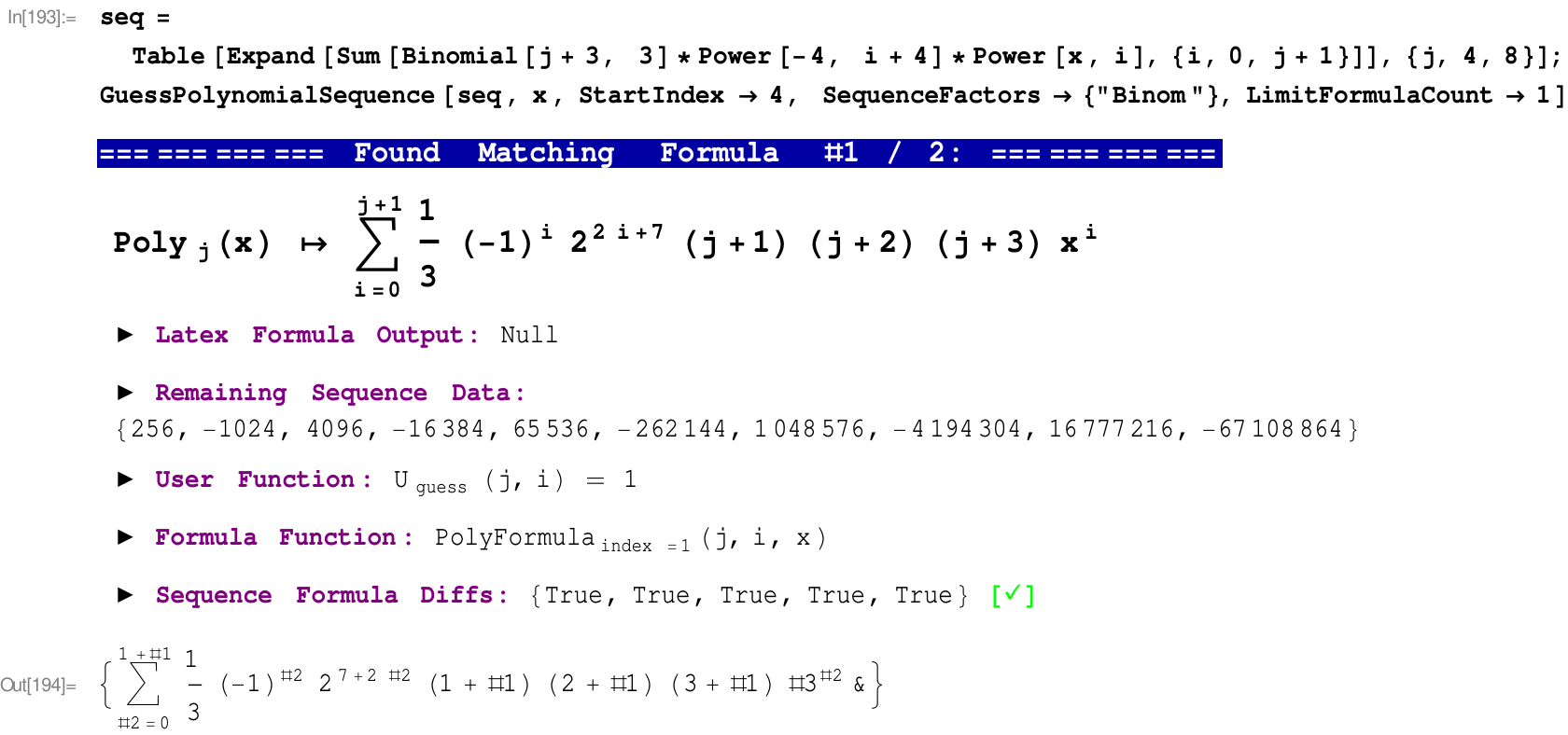}} 
     \end{center} 
     \caption{Expected Sequence Factors Independent of the Sum Index} 
     \label{figure_other-examples-v1} 
\end{figure} 

\begin{figure}[h!] 
     \begin{center}
     \fbox{\includegraphics[width=0.95\textwidth]{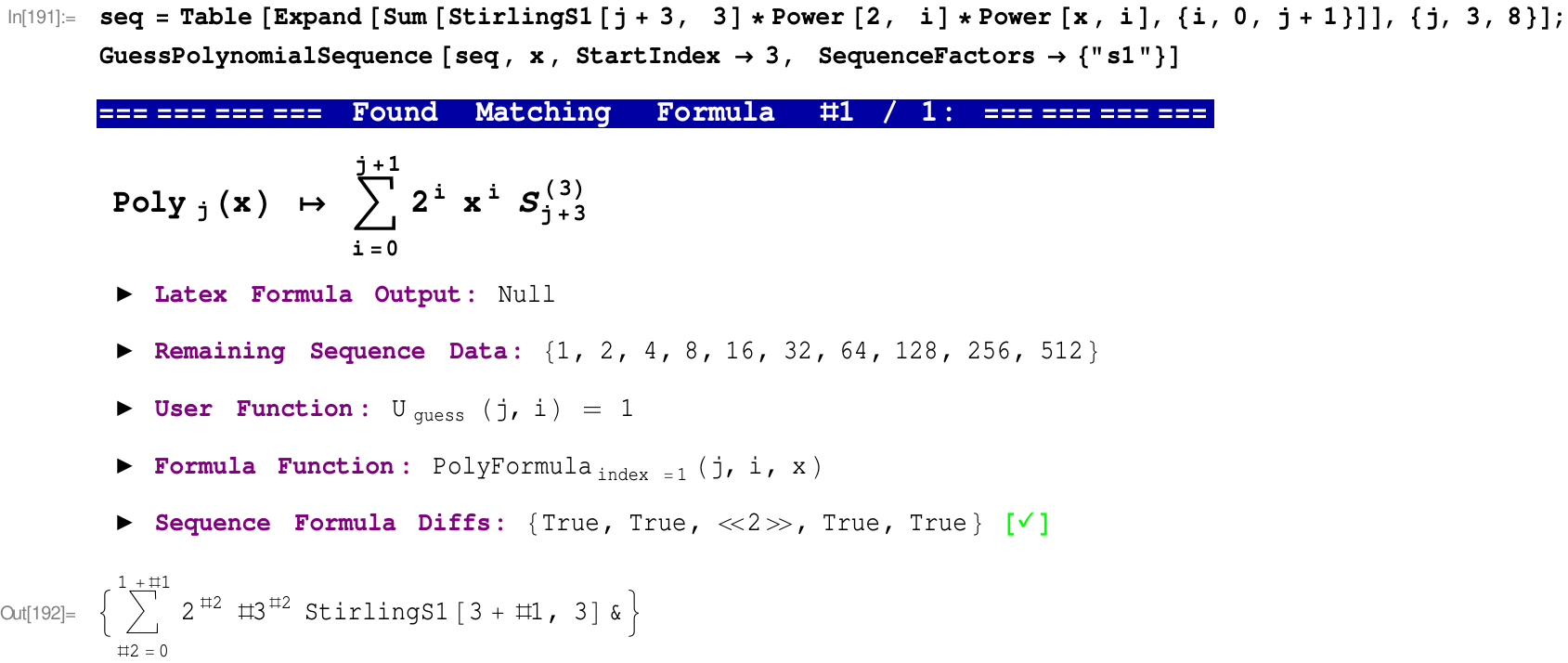}} 
     \end{center} 
     \caption{Another Example of Expected Sequence Factors 
              Independent of the Sum Index} 
     \label{figure_other-examples-v2} 
\end{figure} 

%\subsection{Listings of Other Formulas Obtained by the Package} 
\subsection{Examples: Recognition of Other Sequence Formulas with the 
            \Mm Package} 

Figure \ref{figure_other-examples-v1} and 
Figure \ref{figure_other-examples-v2} 
cite two additional examples of sequence formulas that the package is able to 
recognize when the triangular sequence factors expected by the user do not 
depend on the summation index, $i$, only the polynomial sequence index, $j$. 
In the first example given in 
Figure \ref{figure_other-examples-v1}, the expected binomial coefficient 
factor corresponds to a polynomial in $j$. 
In the second example given in 
Figure \ref{figure_other-examples-v2}, the expected Stirling number 
factor corresponds to an expansion in terms of $r$--order harmonic numbers, 
$H_{j+2}^{(r)}$, that is reported as the factor of the original 
Stirling number sequence. 

\subsection{Examples: Recogntion of Other Polynomial Sum Identities with the 
            \emph{Sage} Package Implementation} 

\subsubsection{Example: A Legendre Polynomial Identity Involving 
               Squares of the Binomial Coefficients} 

A finite polynomial sum over the squared powers of the binomial 
coefficients is expressed through the Legendre polynomials, $P_n(x)$, and 
its ordinary generating function in two variables in the following forms 
\citep[\S 18]{NISTHB}: 
\begin{align*} 
S_n(z) & := \sum_{k=0}^{n} \binom{n}{k}^2 z^k = 
     (1-z)^n P_n\left(\frac{1+z}{1-z}\right) = 
     [t^n] \left(\frac{1}{\sqrt{1 - (1+z) t + (z-1)^2 t^2}}\right). 
\end{align*} 
Alternately, we may obtain information about a closed--form sum for the 
Legendre polynomials over these polynomial inputs to the sequence 
through a known recurrence relation for the sums, $S_n(z)$, given by 
\cite[p. 543]{GKP} 
\begin{equation*} 
(n+1) (1-z)^2 S_n(z) - (2n+3)(z+1) S_{n+1}(z) + (n+2) S_{n+2}(z) = 0. 
\end{equation*} 
Figure \ref{figure_LegendrePoly_SumFormula_formula_v1.2} 
provides a listing of \emph{Sage} commands using the new package 
implementation in \emph{Python} to obtain a sequence formula for the 
right--hand--side polynomial expansions. 

\begin{figure}[h] 
     \begin{center}
     \begin{boxedminipage}{0.87\linewidth}
     \begin{sagecommandline}     
     sage: ## A Legendre polynomial identity involving squared 
     sage: ## powers of the binomial coefficients
     sage: from GuessPolynomialSequenceFunction import * 
     sage: n, z = var('n z')
     sage: poly_seq_func = lambda m: expand( factor( legendre_P(m, (1+z)/(1-z)) * ((1-z) ** m) ) )
     sage: pseq_data = map(poly_seq_func, range(1, 6))
     sage: guess_polynomial_sequence(pseq_data, z, seq_factors = ["Binom2"]);
     \end{sagecommandline}
     \end{boxedminipage}
     \end{center} 
     \caption{An Exponential Generating Function for the 
              Binomial Coefficients (\emph{Sage})} 
     \label{figure_LegendrePoly_SumFormula_formula_v1.2} 
\end{figure} 

\subsubsection{Example: An Exponential Generating Function for the 
               Exponential Bell Polynomials} 

An exponential generating function for the Bell, or exponential, 
polynomials and the corresponding finite sum expansion over the 
Stirling numbers of the second kind is given in the next equation 
\citep[\S 4.1.3]{UC}. 
\begin{align*}
n! \cdot B_n(x) & = [t^n] \exp\left(\left(e^t - 1\right) x\right) 
     = \sum_{k=0}^{n} \gkpSII{n}{k} x^k
\end{align*} 
Figure \ref{figure_BellExpPoly_SumFormula_formula_v1.2} 
provides a listing of the \emph{Sage} commands needed to 
recognize the rightmost identity for this special sequence. 

\begin{figure}[h] 
     \begin{center}
     \begin{boxedminipage}{0.87\linewidth}
     \begin{sagecommandline}     
     sage: ## Finite summation formula for the exponential 
     sage: ## Bell polynomials
     sage: from GuessPolynomialSequenceFunction import *
     sage: def series_coefficient_zpow(f, fvar, ncoeff): return f.taylor(fvar, 0, ncoeff) - f.taylor(fvar, 0, ncoeff - 1)
     sage: def series_coefficient(f, fvar, ncoeff): return series_coefficient_zpow(f, fvar, ncoeff).subs_expr(fvar == 1) 
     sage: n, x, t = var('n x t')
     sage: spoly_ogf = exp( (exp(t) - 1) * x)
     sage: poly_seq_func = lambda n: series_coefficient(spoly_ogf, t, n) * factorial(n)
     sage: pseq_data = map(poly_seq_func, range(1, 6))
     sage: guess_polynomial_sequence(pseq_data, x, seq_factors = ["S2"]);
     \end{sagecommandline}
     \end{boxedminipage}
     \end{center} 
     \caption{An Exponential Generating Function for the 
              Binomial Coefficients (\emph{Sage})} 
     \label{figure_BellExpPoly_SumFormula_formula_v1.2} 
\end{figure} 

\chapter{Conclusions} 
\label{Chapter_Conclusions} 

\section{Concluding Remarks} 

The package source code portion of the thesis 
provides a successful ``proof of concept'' implementation of the logic 
employed by the approach to the package to recognize 
polynomial sequence formula types of the noted forms in 
\eqref{eqn_Polyjx_single_factor_seq_formula_template_v1} and 
\eqref{eqn_Polyjx_double_factor_seq_formula_template_v2}. 
The primary deficiency of the package implementation is current as of 
this writing is the long running time of the package function 
\GPSFPkgGuessFnName when processing double--factor and 
multiple--factor sequence formulas of the form outlined in 
\eqref{eqn_pmx_gen_poly_sequence_forms-intro_stmt_v1}. 
Single--factor polynomial sequence formulas in the form of 
\eqref{eqn_Polyjx_single_factor_seq_formula_template_v1} like those cited in 
\eqref{eqn_ModPolyPowGFs_of_the_S1StirlingNumbers-stmts_v1} of the 
introduction are already somewhat easy, though not trivial, 
to guess by the user. 
For the package to be really useful in practice, the sequence 
recognition routines provided through the wrapper function \GPSFPkgGuessFnName 
should be able to guess double--factor formulas of the form in 
\eqref{eqn_Polyjx_double_factor_seq_formula_template_v2} 
fairly quickly and efficiently out--of--the--box. 

The examples given in \chref{Chapter_PkgUsageChapter} 
provide several non--trivial uses of the package 
for recognizing single--factor polynomial formulas of the 
first sequence form in 
\eqref{eqn_Polyjx_single_factor_seq_formula_template_v1}. 
These and related applications corresponding to polynomials that 
satisfy a single--factor formula of this variety are easily and fairly 
quickly recognized by the package given an accurate user--defined 
setting of the \url{SequenceFactors} runtime option to 
\GPSFPkgGuessFnNamePlain. 

For polynomial sequences that satisfy a double--factor formula of the 
second form in \eqref{eqn_Polyjx_double_factor_seq_formula_template_v2}, and 
more generally a multiple--factor formula in the form stated in 
\eqref{eqn_pmx_gen_poly_sequence_forms-intro_stmt_v1} where $r \geq 3$, the 
current package implementation is unable to quickly search for matching 
formulas without a somewhat manual limited setting of the 
\url{IndexOffsetPairs} option provided at runtime. 
The sample output for the examples given in 
Figure \ref{figure_S1S2-double-factor-formula-example-v1} and 
Figure \ref{figure_S1Binom-double-factor-formula-example-v2} 
show the usage of the package for handling double--factor sequence 
formulas with an appropriate setting of this option. 
%However, 
In future revisions of the package, 
it should ideally be possible for the package to quickly obtain 
formulas for these sequence cases without the user manually resetting the 
default search options used with the \GPSFPkgGuessFnName 
function provided by the package. 

\section{Future Features in the Package} 
\label{Chapter_PkgImpl_Section_FutureFeatures} 

\subsection{Processing Polynomial Sequences with Rational Coefficients} 
\label{subSection_FutureFeaturesChapter_RationalSeqTransforms} 

One approach to extending the package functionality to recognize 
formulas for polynomial sequences in $\mathbb{Q}\lbrack x\rbrack$ 
is to pre--process the rational--valued coefficients to 
transform the sequence into the polynomials over the integers already 
handled by the package routines. Variations of these pre--processing 
transformations include normalizing the polynomials, its coefficients, or 
both by exponential factors to clear the denominators of the 
rational--valued input sequence. For example, let the polynomial 
$p_j(x) := \sum_i c_i x^i$. Then these transformations are formulated as 
obtaining the modified polynomials, $\widetilde{p}_j(x)$, as 
$\widetilde{p}_j(x) := j! \cdot p_j(x)$, as 
$\widetilde{p}_j(x) := \sum_i i! \cdot c_i x^i$, or in the combined form of 
$\widetilde{p}_j(x) := \sum_i j! \cdot i! \cdot c_i x^i$, 
whenever the resulting modified polynomial sequences are in 
$\mathbb{Z}\lbrack x \rbrack$. 

Another transformation option is applied to rationalize the polynomial 
sequences, $S_m(n)$, in $n$ defined through the following sums where 
$B_n$ denotes the (rational) sequence of Bernoulli numbers 
\citep[\S 6.5]{GKP}: 
\begin{equation} 
\notag 
S_m(n) := \sum_{k=0}^{n-1} k^m = \sum_{k=0}^{m} \binom{m+1}{m-k} 
     \frac{B_{m-k}}{(m+1)} \cdot n^{k+1}. 
\end{equation} 
The sequences in the previous equation are normalized by multiplying 
each polynomial, $S_m(n)$, by the least common multiple of the 
denominators of each coefficient of $n^{k+1}$ in the formula. 
Then assuming access to the lookup capabilities of the \OEIS database,
which contains sequence entries for both integer sequences of the 
numerators and denominators of the Bernoulli numbers, obvious 
factors, of say $691$, are recognized to process the full formulas for the 
sequences of $S_m(n)$ over $\mathbb{Q}\lbrack n \rbrack$. 

\subsection{Polynomial Expansions With Respect to a Suitable Basis} 

The discussion given in \citep[Appendix A]{RATE.M-PKG-DOCS} related to the 
implementation of the \ttemph{Rate} package for 
\Mm states a useful observation that may be 
adapted to the polynomial formula searches local to this package. 
Specifically, expressing input polynomial sequences with respect to a 
``\emph{suitable}'' basis, like shifted factorial functions or 
polynomial terms expressed by binomial coefficients, allows for 
recognition of sequence formulas that are not apparent in the default 
expansions of the polynomial sequence variable. 
Several examples relevant to adapting this idea in the context of the 
factorization--based approach in this package include the following 
polynomial sequence expansions \citep[Ex. 6.78; \S 6.2; Ex. 6.68]{GKP}: 
\begin{align*} 
\binom{2n}{n} \frac{B_n}{(n+1)} & = \sum_{k=0}^{n} \gkpSII{n+k}{k} 
     \binom{2n}{n+k} \frac{(-1)^{k}}{(k+1)} \\ 
x^n & = \sum_{k=0}^{n} \gkpEI{n}{k} \binom{x + k}{n} \\ 
\gkpSI{x}{x-n} & = \sum_{k \geq 0} \gkpSII{n}{k} \binom{x+k}{2n},\ n \geq 0 \\ 
\gkpEII{n}{m} & = \sum_{k=0}^{m} 
     \binom{2n+1}{k} \gkpSII{n+m+1-k}{m+1-k} (-1)^{k},\ n > m \geq 0. 
\end{align*} 
These sequences provide applications related to the 
polynomial expansions of the Catalan numbers (in $n$), the 
Stirling convolution polynomials, $\sigma_n(x)$, and the 
second--order Eulerian numbers, $\gkpEII{n}{m}$, respectively. 

%\section{Topics for Future Research} 
\section{Future Research Topics} 

The next sections discuss several topics for future research 
suggested by the implementation of the software package for the thesis. 
These future research topics include 
a new variation of integer factorization algorithms 
motivated by the factorization--based approach to handling the 
user--defined expected sequence factors in the package routines, as well as 
additional topics for future exploration to extend the current 
capabilities of the univariate polynomial sequence recognition in the 
package. 
%% 
%Notice that 
The extension of the current package functionality to 
recognizing polynomials in a single variable with rational--valued 
coefficients is already considered in 
\sref{Chapter_PkgImpl_Section_FutureFeatures} of the thesis above. 
%Potential 
%Optimizations of the current source code for the package is also 
%already discussed in previous sections of the thesis. 

%\section{New Sequence--Based Factorization Algorithms Suggested by the Package} 
\subsection{Sequence--Based Integer Factorization Algorithms} 

The treatment of the user--defined expected sequence factors as 
``primitives'' in the formulas returned by the package functions 
motivates the construction of a class of integer factorization algorithms 
formulated briefly in the discussion below. 
Much like computing the prime factorization of an arbitrary integer, 
this class of algorithms should compute the decomposition of an 
integer into a product of elements over some specified set of integer 
sequences where the elements of these sequences are treated as ``atoms'' in 
the factorization returned by the procedure. 

%% 
%For example, 
Stated more precisely: 
given an integer $i$ (or some set of integer--valued 
polynomial coefficients) and a list of $k$ integer sequences, 
$\{S_1, S_2, \ldots, S_k\}$, 
we seek the most efficient way to decompose the integer into %a list of 
all possible products of integer factors of the form 
\begin{equation} 
%\notag 
\label{eqn_IntFactorAlg_int_decomp_form-stmt_v1} 
i := f_1 \cdot f_2 \cdot \cdots \cdot f_k \times r, 
\end{equation} 
where the factor $f_i$ belongs to the sequence $S_i$ 
(for each $1 \leq i \leq k$), and where the remaining factor term, $r$, 
is reserved for later processing. 
The computation of the list of all factors of the form in 
\eqref{eqn_IntFactorAlg_int_decomp_form-stmt_v1} can be computed 
over some specified number of elements of each sequence, 
or a fixed number of rows for the case of a triangular sequence, $S_i$. 
%especially if the sequences involved contain infinitely--many $1$'s or 
It seems reasonable to expect that such an algorithms must employ the 
prime factorizations of the individual factor sequences, $F_i$. 
We also seek a solution in the general case, though of course it may be 
possible to derive sequence--specific procedures, say, to recognize 
factors of the Stirling number or binomial coefficient triangles. 

The need for this type of factorization is apparently new, as searches 
for such subroutines to employ within the package returned no useful 
known results, though 
it is possible that there are existing prime factorization algorithms 
that may be especially well--suited, or adapted, to this purpose. 
This required factorization procedure is handled as an 
inefficient implementation of an oracle of sorts within the current 
implementation of the \Mm package. 

\bibliographystyle{abbrv-mod}

\begin{thebibliography}{10}

\bibitem{ACOMB-BOOK}
P.~Flajolet and R.~Sedgewick, {\em Analytic Combinatorics}, Cambridge
  University Press, 2009.

\bibitem{GKP}
R.~L. Graham, D.~E. Knuth, and O.~Patashnik, {\em Concrete Mathematics: A
  Foundation for Computer Science}, Addison-Wesley, 1994.

\bibitem{AXIOM-GUESS-PKG-DOCS}
W.~Hebisch and M.~Rubey, Extended rate, more {G}{F}{U}{N}, {\em Journal of
  Symbolic Computation} {\bf 46} (2011),  889--903.
\newblock See
  \url{http://axiom-wiki.newsynthesis.org/GuessingFormulasForSequences}.

\bibitem{STIRLING.M-PKG-DOCS}
M.~Kauers, Summation algorithms for {S}tirling number identities, {\em Journal
  of Symbolic Computation} {\bf 42} (2007),  948--970.
\newblock See
  \url{http://www.risc.jku.at/research/combinat/software/ergosum/RISC/Stirling.html}.

\bibitem{RATE.M-PKG-DOCS}
C.~Krattenthaler, Advanced determinant calculus, {\em S{\'e}minaire
  Lotharingien de Combinatoire} {\bf 42} (1999),  1--67.
\newblock See \url{http://www.mat.univie.ac.at/~kratt/rate/rate.html}.

\bibitem{NISTHB}
F.~W.~J. Olver, D.~W. Lozier, R.~F. Boisvert, and C.~W. Clark, eds., {\em
  {NIST} Handbook of Mathematical Functions}, Cambridge University Press, 2010.

\bibitem{AEQUALSB-BOOK}
M.~Petkovsek, H.~S. Wilf, and D.~Zeilberger, {\em {A} = {B}}, {A} {K} {P}eters,
  Ltd., 1996.

\bibitem{UC}
S.~Roman, {\em The Umbral Calculus}, Dover, 1984.

\bibitem{GFUN-PKG-DOCS}
B.~Salvy and P.~Zimmermann, Gfun: a {M}aple package for the manipulation of
  generating and holonomic functions in one variable, {\em ACM Transactions on
  Mathematical Software} {\bf 20} (1994),  163--177.

\bibitem{OEIS-UPDATED-URL}
N.~J.~A. Sloane, The {O}nline {E}ncyclopedia of {I}nteger {S}equences, \\ {\tt
  http://oeis.org}.

\end{thebibliography}

\end{document}